\documentclass{amsart}
\usepackage{color}
\usepackage{graphicx,epsf}
\newtheorem{theorem}{Theorem}[section]
\newtheorem{lemma}[theorem]{Lemma}

\theoremstyle{definition}
\newtheorem{definition}[theorem]{Definition}

\theoremstyle{remark}

\numberwithin{equation}{section}

\newtheorem{proposition}{Proposition}
\newcommand {\sign} {\mbox{\bf sign}}

\newtheorem{algorithm}{\sc Algorithm}

\newcommand {\R} {\mbox{\rm R}}
\newcommand {\eop}      {\hfill $\Box$}
\newcommand {\junk}[1]{}
\newcommand {\hide}[1]{}

\newcommand {\s}        {\mbox{\rm sign}}
\newcommand {\D}     {\mbox{\rm D}}

\newcommand {\Real}[1]   {\mbox{${\mathbb R}^{#1}$}}

 \newcommand {\re}         {\Real{}}

\newcommand {\Sign}         {\mbox{\rm SIGN}}

\newcommand {\Z}  {{\mathbb Z}}
 
\newcommand {\Q}         {{\mathbb Q}}

\newcommand {\ZZ} {{\rm Z}}
\newcommand {\RR} {{\mathcal R}}

\newcommand {\la}   {{\langle}}
\newcommand {\ra}   {{\rangle}}
\newcommand {\eps} {{\varepsilon}}
\newcommand {\E} {{\rm Ext}}

\def\sign{{\rm sign}}

\newcommand {\Ker}      {\mbox{\rm Ker}}
\newcommand {\Ima}      {\mbox{\rm Im}}

\def\addots{\mathinner{\mkern1mu
\raise1pt\vbox{\kern7pt\hbox{.}}
\mkern2mu\raise4pt\hbox{.}\mkern2mu
\raise7pt\hbox{.}\mkern1mu}}

\newcommand{\RM}  {\mbox{\rm RM}}

\renewcommand {\Im} {\mbox {\rm Im}}

\newcommand{\coucou}[1]{\ifvmode\else\marginpar[\hfill$\rhd$]{$\lhd$}\fi
                       $\langle$\textsc{#1}$\rangle$}

\newcommand{\basu}[1]{}
\newcommand {\rp}[1]{}

\begin{document}

\title[Computing the First Betti Number]{Computing 
the First Betti Number
and Describing the Connected Components
of Semi-algebraic Sets
\footnote{2000 Mathematics Subject Classification 14P10, 14P25}}

\author{Saugata Basu}
\address{School of Mathematics,
Georgia Institute of Technology, Atlanta, GA 30332, U.S.A.}
\email{saugata@math.gatech.edu}
\thanks{Supported in part by an NSF Career Award 0133597 and
a Alfred P. Sloan Foundation Fellowship.}

\author{Richard Pollack}
\address{Courant Institute of Mathematical Sciences, 
New York University, New York, NY 10012, U.S.A.}
\email{pollack@cims.nyu.edu}
\thanks{Supported in part by 
NSA  grant MDA904-01-1-0057 and NSF grants
CCR-9732101 and  CCR-0098246.}

\author{Marie-Fran\c{c}oise Roy}
\address{IRMAR (URA CNRS 305), 
Universit\'{e} de Rennes,
Campus de Beaulieu 35042 Rennes cedex FRANCE.}
\email{marie-francoise.roy@univ-rennes1.fr}



\keywords{Betti numbers, Semi-algebraic sets}

\begin{abstract}
In this paper we  describe a singly exponential algorithm for
computing the first Betti number of a given 
semi-algebraic set.   Singly
exponential algorithms for computing the zero-th Betti number, and
the Euler-Poincar\'e characteristic, were known before. No singly exponential
algorithm was known for computing any of the individual Betti numbers
other than the zero-th one.
We also give
algorithms for obtaining  semi-algebraic descriptions
of the semi-algebraically connected components of any given real algebraic
or semi-algebraic set in single-exponential time improving on previous
results.
\end{abstract}

\maketitle
\section{Introduction}
\label{sec:intro}
Let $\R$ be a real closed field and $S \subset \R^k$ a semi-algebraic set
defined by a
quantifier-free
Boolean formula with atoms of the form
$P > 0, P < 0, P=0$ for $P \in {\mathcal P} \subset \R[X_1,\ldots,X_k]$.
We call $S$ a ${\mathcal P}$-semi-algebraic set.
If, instead, the Boolean formula
has atoms of the form $P=0, P \geq 0, P \leq 0, \;P \in {\mathcal P}$, 
and additionally contains no negation,
then we will call $S$ a ${\mathcal P}$-closed semi-algebraic set.
It is well known \cite{O,OP,Milnor,Thom,B99,GV} 
that the topological complexity of $S$ 
(measured by the various Betti numbers of $S$) is bounded by $O(sd)^k$,
where $s  = \#({\mathcal P})$ and $d = \max_{P\in {\mathcal P}}{\rm deg}(P).$
More precise bounds on the individual 
Betti numbers of $S$ appear in \cite{B03}.
Even though the Betti numbers of $S$ are bounded singly exponentially
in $k$, 
there is no singly exponential algorithm for computing the Betti numbers of
$S$. 
This absence is related to the fact that
there is no known algorithm for producing a singly exponential
sized triangulation of $S$ (which would immediately imply a singly exponential
algorithm for computing the Betti numbers of $S$). In fact, the
existence of a singly exponential sized triangulation, is considered to
be a major open question in 
algorithmic
real algebraic geometry. 
Moreover, determining the exact complexity of computing the Betti numbers
of semi-algebraic sets is an area of active research in computational 
complexity theory, for instance counting versions
of complexity classes in the Blum-Shub-Smale model of computation
(see \cite{BC}).

Doubly exponential algorithms (with complexity $(sd)^{2^{O(k)}}$)
for computing all the Betti numbers are known, 
since it is possible to obtain a triangulation of $S$ in doubly
exponential time using cylindrical algebraic decomposition 
\cite{Col,BPR03}.
In the absence of a singly exponential time algorithm for computing
triangulations of semi-algebraic sets, algorithms with single exponential
complexity are known only for the problems of testing emptiness
\cite{R92,BPR95},
computing the zero-th Betti number 
(i.e. the number of semi-algebraically connected components of 
$S$) \cite{GV92,Canny93a,GR92,BPR99},
as well as the Euler-Poincar\'e characteristic of $S$ \cite{B99}.

In this paper we describe the first singly exponential algorithm for
computing the first Betti number of a given 
semi-algebraic set  $S \subset \R^k$. 
In the process,
we also give efficient 
algorithms for obtaining semi-algebraic descriptions 
of the semi-algebraically connected components of a given real algebraic or 
semi-algebraic set.
These algorithms have complexity bounds which improve the
complexity of the best previously known algorithm \cite{HRS94}.

The rest of the paper is organized as follows.

There are several ideas involved in the design of our algorithm for computing the
first Betti number of a given semi-algebraic set, which corresponds to
the  main steps in our algorithm. We describe each of them separately in 
different sections.

In Section \ref{15:sec:desc} we recall the notion of a roadmap
of an algebraic set \cite{BPR99}
and indicate how to use it to construct connecting paths in basic
semi-algebraic sets.

In Section 
\ref{sec:covering}  we define
certain semi-algebraic sets which we call parametrized paths and
prove that 
under a certain hypothesis these sets  are semi-algebraically contractible.
We also outline the input, output, and complexity
of an algorithm computing 
a covering of a given basic semi-algebraic set, $S \subset \R^k$,
by a singly exponential number of 
parametrized paths. 

In Section \ref{sec:acycliccov}, we use the properties of parametrized
paths proved in Section \ref{sec:covering} to give an algorithm 
(Algorithm \ref{16:alg:acycliccovering})
for computing a covering of a given closed and bounded semi-algebraic set 
by a single exponential sized family of closed, bounded as well as
contractible semi-algebraic sets. The complexity of this algorithm is
singly exponential.

In Section \ref{sec:top}, we recall some results from algebraic
topology which allows us 
to compute the first Betti number of a closed and bounded
semi-algebraic set from a covering 
of the given set consisting of
closed, bounded and contractible sets.
The main tool 
here is a spectral sequence associated to the
Mayer-Vietoris double complex. We show how to compute the
first Betti number once we have computed a covering by
closed contractible sets and the number
of connected components of their pair-wise and triple-wise intersections of 
the  sets in this covering and their incidences. 
If the size of the  covering is singly exponential,
this yields a singly exponential  algorithm for computing the first 
Betti number. Extensions of these ideas for computing 
a fixed number of higher Betti numbers 
in singly exponential time is 
possible and  is reported on in a subsequent paper \cite{B04}.
In Section \ref{sec:closedcase}, we describe an algorithm for computing the
first Betti number of a given ${\mathcal P}$-closed semi-algebraic set.

In Section \ref{sec:GV} we recall a technique introduced 
by Gabrielov and Vorobjov \cite{GV},
for replacing any given semi-algebraic set by one which
is closed and bounded and has the same homotopy type as the given set. 
In fact, we prove a slight strengthening of the
main result in \cite{GV}, in that we prove that 
the new set has the same homotopy type as the given one, while
the 
corresponding result (Lemma 5) 
in \cite{GV} states that only the sum of the Betti numbers
is preserved.
The above construction  allows us
to reduce the case of general semi-algebraic sets to ones
which are closed and bounded treated in Section \ref{sec:top}
without any significant worsening of complexity.

In Section \ref{sec:general}, we describe an algorithm for computing the 
first Betti number of a general semi-algebraic set, using the construction
described in Section \ref{sec:GV} to first reduce the problem to the
${\mathcal P}$-closed case already treated in Section \ref{sec:closedcase}.

Finally, in Section \ref{sec:cc} we indicate that the 
algorithms described in Section \ref{sec:covering} 
actually produces descriptions of the connected
components of a given algebraic or semi-algebraic set
in an efficient manner. 

\section{Preliminaries}
\label{sec:prelim}
Let $\R$ be  a real closed field. For an element $a \in\R$ we let
$$
\s(a) =
\begin{cases}
0 & \mbox{ if }  a=0,\cr
1 & \mbox{ if } a> 0,\cr
-1& \mbox{ if } a< 0.
\end{cases}
$$
If ${\mathcal P}$ is a finite subset of
$\R [X_1, \ldots , X_k]$, we write the {\em set of zeros}
of ${\mathcal P}$ in $\R^k$ as
$$
\ZZ({\mathcal P},\R^k)=\{x\in \R^k\mid\bigwedge_{P\in{\mathcal P}}P(x)= 0\}.
$$
We denote by 
$B(0,r)$ the open ball with center 0 and radius $r$.

Let ${\mathcal Q}$ and ${\mathcal P}$ be finite subsets of 
$\R[X_1,\ldots,X_k],$ $Z = \ZZ({\mathcal Q},\R^k),$ and
$Z_r = Z \cap B(0,r).$
A  {\em sign condition}  on
${\mathcal P}$ is an element of $\{0,1,- 1\}^{\mathcal P}$.
The {\em realization of the sign condition
$\sigma$ over $Z$}, $\RR(\sigma,Z)$, is the basic semi-algebraic set
\label{def:R(Z)}
$$
        \{x\in \R^k\;\mid\; \bigwedge_{Q \in {\mathcal Q}} Q(x)=0
     \wedge \bigwedge_{P\in{\mathcal P}} \s({P}(x))=\sigma(P) \}.
$$
The {\em realization of the sign condition
$\sigma$ over $Z_r$}, $\RR(\sigma,Z_r),$ is the basic semi-algebraic set
\label{def:R(Z_r)}
$
        \RR(\sigma,Z) \cap B(0,r).
$
For the rest of the paper, we fix an open
ball $B(0,r)$ with center 0 and radius $r$ big enough
so that, for  every sign condition $\sigma$,
$\RR(\sigma,Z)$ and $\RR(\sigma,Z_r)$ are homeomorphic.
This is always possible by the local conical structure at infinity
of semi-algebraic sets 
(\cite{BCR}, page 225).

A closed and bounded semi-algebraic
set $S \subset \R^k$ is semi-algebraically triangulable (see \cite{BPR03}),
and we denote by $H_i(S)$  the $i$-th simplicial homology group of $S$
with rational coefficients.  The groups  $H_i(S)$ are invariant under
semi-algebraic homeomorphisms and coincide with the corresponding
singular homology groups when $\R = \re$. We denote by $b_i(S)$ the
$i$-th Betti number of $S$ (that is, the dimension  of $H_i(S)$ as a vector
space), and by $b(S)$ the sum $\sum_i b_i(S)$.
For a closed but not necessarily bounded semi-algebraic set $S \subset \R^k$,
we will denote by $H_i(S)$  the $i$-th simplicial homology group of 
$S \cap \overline{B(0,r)}$, where $r$ is sufficiently large. 
The sets $S \cap \overline{B(0,r)}$ are semi-algebraically homeomorphic
for all sufficiently large $r> 0$, by
the local conical structure at infinity of 
semi-algebraic sets, 
and hence this definition makes sense.

The definition of homology groups of arbitrary semi-algebraic sets in
$\R^k$ requires some care and several possibilities exist. In this paper,
we define the homology groups of realizations of sign conditions as follows.

Let $\R$ denote a real closed field and 
$\R'$ a real closed field containing $\R$.
Given a semi-algebraic set
$S$ in ${\R}^k$, the {\em extension}
of $S$ to $\R'$, denoted $\E(S,\R'),$ is
the semi-algebraic subset of ${ \R'}^k$ defined by the same
quantifier free formula that defines $S$.
The set $\E(S,\R')$ is well defined (i.e. it only depends on the set
$S$ and not on the quantifier free formula chosen to describe it).
This is an easy consequence of the transfer principle \cite{BPR03}.

Now, let $S \subset \R^k$ be a ${\mathcal P}$-semi-algebraic set, where
${\mathcal P} = \{P_1,\ldots,P_s \}$ is a finite subset of $\R[X_1,\ldots,X_k].$ 
Let $\phi(X)$ be a quantifier-free formula defining $S$. 
Let $P_i = \sum_{\alpha} a_{i,\alpha}X^\alpha$ where the $a_{i,\alpha} \in \R.$
Let
$A = (\ldots,A_{i,\alpha},\ldots )$ denote the vector of variables 
corresponding 
to the coefficients of  the polynomials in the family ${\mathcal P},$
and let
$a = (\ldots,a_{i,\alpha},\ldots) \in \R^N$ denote the vector of
the actual coefficients of the polynomials in ${\mathcal P}$. 
Let $\psi(A,X)$ denote the formula obtained
from $\phi(X)$ by replacing each coefficient of each polynomial in ${\mathcal P}$
by the corresponding variable, so that $\phi(X) = \psi(a,X).$ It follows from
Hardt's triviality theorem for semi-algebraic mappings \cite{Hardt},
that there exists,
$a' \in \re_{\rm alg}^N$ such that 
denoting by $S' \subset \re_{\rm alg}^k$ the semi-algebraic set defined by
$\psi(a',X)$, the semi-algebraic set
$\E(S',\R)$, has the same homeomorphism type as $S$. We define the 
homology groups of $S$ to be the singular homology groups 
of $\E(S',\re).$ It follows from the Tarski-Seidenberg transfer principle,
and the corresponding property of singular homology groups,
that the homology groups defined this way are invariant under 
semi-algebraic homotopies. It is also clear that this
definition is  compatible with the
simplicial homology for closed, bounded semi-algebraic sets, and 
the singular homology groups when the ground field is $\re$.
Finally it is also clear that, the Betti numbers are not changed after 
extension:
$b_i(S)=b_i(\E(S,\R')).$

\section{Roadmap of a semi-algebraic set}
\label{15:sec:desc}
We first define a roadmap 
of a semi-algebraic set. Roadmaps are crucial ingredients in all singly
exponential algorithms known for computing connectivity properties
of semi-algebraic sets
such as computing the number of connected 
components, as well as testing whether two points of a given semi-algebraic 
set belong to the same semi-algebraically connected component.

We use the following notations.
Given  $x=(x_1,\ldots,x_k)$ we write  $\bar x_i$ for $(x_1,\ldots,x_i),$
and  $\tilde x_i$ for $(x_{i+1},\ldots,x_k).$
We also denote by $\pi_{1\ldots j}$ the projection,
$x \mapsto \bar x_j.$
Given a set $S \subset \R^k$, $y \in \R^j$ we denote by 
$S_y=S \cap \pi_{1 \ldots j}^{-1}(y)$.

\label{15:def:roadmap}
Let $S \subset \R^k$ be a semi-algebraic set.
A {\em roadmap}
\index{Roadmap} for $S$
 is a semi-algebraic set $M$
of dimension at most one contained in $S$
which satisfies the following roadmap
conditions:
\begin{itemize}
\item ${\rm RM}_1$ For every semi-algebraically
connected component
$D$ of $S$,
$D \cap M$ is semi-algebraically connected.
\item ${\rm RM}_2$ For every $x \in {\R}$ and
for every semi-algebraically connected component $D'$
of $S_x$, $D'\cap  M \neq \emptyset.$
\end{itemize}

We describe the construction of a roadmap
$M$ for a bounded algebraic set $\ZZ(Q,\R^k)$
 which contains a finite set of
points ${\mathcal N}$ of $\ZZ(Q,\R^k)$. A precise description of how the
construction can be performed algorithmically
can be found in \cite{BPR03}.

A key ingredient of the roadmap is  the construction of 
a particular finite set of points having the property that, 
they intersect every connected component of $\ZZ(Q,\R^k)$.
We call them $X_1$-pseudo-critical points,
since they are obtained as limits of the  critical points of the projection to
the $X_1$ coordinate of  a bounded nonsingular
algebraic hypersurface defined by a particular
infinitesimal
deformation
of the polynomial $Q$.
Their projections on the $X_1$-axis are called pseudo-critical values.
These points are obtained as follows. 

We denote by $\R\langle \zeta\rangle$  the real closed field of algebraic
Puiseux series in $\zeta$ with coefficients in $\R$ \cite{BPR03}. 
The sign of a Puiseux series in $\R\langle \zeta\rangle$
agrees with the sign of the coefficient
of the lowest degree term in
$\zeta$. 
This induces a unique order on $\R\langle \zeta\rangle$ which
makes $\zeta$
infinitesimal: $\zeta$ is positive and smaller than
any positive element of $\R$.
When $a \in \R\la \zeta \ra$ is bounded by an element of $\R$,
$\lim_\zeta(a)$ is the constant term of $a$, obtained by
substituting 0 for $\zeta$ in $a$.
We now define the 
deformation $\bar Q$ of $Q$ as follows.
Suppose that $\ZZ(Q,\R^k)$ is contained in the ball of center
$0$ and radius $1/c$.
 Let
$\bar d$ be an even integer bigger than the degree $d$ of $Q$,
\begin{equation}
G_k(\bar d,c)
=c^{\bar d}(X_1^{\bar d}+\cdots+
X_{k}^{\bar d}+X_2^2+\cdots+X _k^2)-(2k-1),
\end{equation}
\begin{equation}
\label{11:equation:deform1}
\bar Q=\zeta G_k(\bar d,c)+{(1-\zeta)
}Q.
\end{equation}

The algebraic
set $\ZZ(\bar Q,\R \la \zeta \ra ^k)$
is a bounded and non-singular
hypersurface lying infinitesimally close to $\ZZ(Q,\R^k)$, 
and the critical points of the projection
map onto the $X_1$ co-ordinate restricted to
$\ZZ(\bar Q,\R \la \zeta \ra ^k)$
form a finite set of points.
We take the 
images of these points under $\lim_\zeta$ and
we call the points obtained in this manner  the $X_1$-pseudo-critical points
of $\ZZ(Q,\R^k)$.
Their projections on the $X_1$-axis  are called pseudo-critical values.

The construction of the roadmap of an algebraic set
containing a finite number 
of input points ${\mathcal N}$ of this algebraic set
is as follows.
We first construct
$X_2$-pseudo-critical points on $\ZZ(Q,\R^k)$
in a parametric way along
 the $X_1$-axis, by following continuously, as $x$ varies on the $X_1$-axis,
the $X_2$-pseudo-critical points on $\ZZ(Q,\R^k)_{x}$. This results in curve segments
and their endpoints on $\ZZ(Q,\R^k).$ The curve segments
 are continuous semi-algebraic curves
parametrized by open intervals on the $X_1$-axis, and their
 endpoints are points of $\ZZ(Q,\R^k)$
above the corresponding endpoints of the open intervals.
Since these curves and their endpoints include,
for every $x\in\R$, the
$X_2-$pseudo-critical points of
$\ZZ(Q,\R^k)_{x}$,
they meet every
connected component of
$\ZZ(Q,\R^k)_{x}$.  Thus the set of curve
segments and their endpoints already satisfy
${\rm RM}_2.$ However, it is clear that this set might not
be semi-algebraically connected in a semi-algebraically
connected component, so ${\rm RM}_1$
might not be
satisfied.
We add additional curve segments to ensure
connectedness by recursing in certain
 distinguished hyperplanes
defined by
$X_1=z$ for
distinguished values
$z.$

The set of {\em distinguished values}
\index{Distinguished!value}
\label{15:def:distinguishedvalues}
is the union
 of the $X_1$-pseudo-critical values, the first
coordinates of the input points ${\mathcal N}$ and the
first coordinates of the endpoints of the curve
segments. A {\em distinguished hyperplane}
\index{Distinguished!hyperplane}
is an hyperplane defined by $X_1=v$, where
$v$ is a distinguished value. The input points, the endpoints of
the curve segments and the intersections of the curve
segments with the distinguished hyperplanes define
the set of {\em distinguished points}.
\label{15:def:distinguishedpoints}.

Let
the distinguished values be
$v_1<\ldots <v_\ell.$
Note that amongst these  are the $X_1$-pseudo-critical values. Above each
interval $(v_i, v_{i+1})$, we have  constructed a collection
 of curve segments ${\mathcal C}_i$ meeting every
semi-algebraically connected component of
$\ZZ(Q,\R^k)_v$ for every $v \in (v_i, v_{i+1})$. Above
each distinguished value $v_i$,  we have 
a set of distinguished points ${\mathcal N}_i$. 
Each curve segment in ${\mathcal C}_i$
has an endpoint in ${\mathcal N}_i$ and another in
${\mathcal N}_{i+1}$.  Moreover, the union of the ${\mathcal N}_i$
contains ${\mathcal N}$.

We then repeat this construction in each 
distinguished hyperplane $H_i$ defined by $X_1=v_i$
with input $Q(v_i,X_2,\ldots,X_k)$ and the
distinguished points in ${\mathcal N}_i$. 
Thus, we construct distinguished values, 
$v_{i,1},\ldots, v_{i,\ell(i)}$ of 
$\ZZ(Q(v_i,X_2,\ldots,X_k),\R^{k-1})$ 
(with the role of $X_1$ being now played by $X_2$) and 
the process is iterated until
for 
$I=(i_1,\ldots,i_{k-2}), 1 \leq i_1 \leq\ell,
\ldots, 1 \leq i_{k-2} \leq \ell(i_1,\ldots,i_{k-3}),$ we have
distinguished values
$v_{I,1}< \ldots < v_{I, \ell(I)}$ along the
$X_{k-1}$ axis with corresponding sets of curve
segments and sets of distinguished points with the
required incidences between them.

The following proposition is proved in \cite{BPR99}
(see also \cite{BPR03}).
\begin{proposition}
\label{15:prop:rm}
The semi-algebraic set
$\RM(\ZZ(Q,\R^k),{\mathcal N})$ obtained by this construction is a
roadmap for $\ZZ(Q,\R^k)$ containing ${\mathcal N}$.
\end{proposition}

Note that if $x \in \ZZ(Q,\R^k)$,  $\RM(\ZZ(Q,\R^k),\{x\})$ 
contains a path,
$\gamma(x)$,
 connecting
a distinguished point $p$ 
of   $\RM(\ZZ(Q,\R^k))$
to $x$.

Later in this paper we shall 
examine the properties of parametrized paths which are the
unions of connecting paths starting 
at a given $p$ and ending 
at $x$,
where $x$ varies over a certain semi-algebraic subset of   
$\ZZ(Q,\R^k)$.  
In order to do so
it is useful to have a better understanding of the
structure of these connecting paths -- especially, of their dependence
on $x$.

Recall that given  $x=(x_1,\ldots,x_k)$ 
we write  $\bar x_i$ for $(x_1,\ldots,x_i),$
and  $\tilde x_i$ for $(x_{i+1},\ldots,x_k).$
We first note that for any $x = (x_1,\ldots,x_k) \in \ZZ(Q,\R^k)$, 
we have by construction that,  $\RM(\ZZ(Q,\R^k))$ is contained in
$\RM(\ZZ(Q,\R^k),\{x\})$. In fact,
$$
\displaylines{
\RM(\ZZ(Q,\R^k),\{x\}) = 
\RM(\ZZ(Q,\R^k)) \cup \RM(\ZZ(Q,\R^k)_{x_1}, {\mathcal M}_{x_1}),
}
$$
where ${\mathcal M}_{x_1}$  
consists of
$\tilde x_1$ 
and the finite set of points obtained by intersecting
the curves in $\RM(\ZZ(Q,\R^k))$ parametrized by the $X_1$-coordinate,
with the hyperplane $\pi_{1}^{-1}(x_1)$.

A connecting path $\gamma(x)$ 
(with non-self intersecting image) 
joining a distinguished point $p$ of $\RM(\ZZ(Q,\R^k))$ 
to $x$ can be extracted from $\RM(\ZZ(Q,\R^k),\{x\})$.
The connecting path $\gamma(x)$ consists of two consecutive parts,  
$\gamma_0(x)$ and $\Gamma_1(x)$.
The  path  $\gamma_0(x)$ is contained in $\RM(\ZZ(Q,\R^k))$
and the path $\Gamma_1(x)$  is contained in $\ZZ(Q,\R^k)_{x_1}$.
The part $\gamma_0(x)$ consists of a sequence of sub-paths,
$\gamma_{0,0},\ldots,\gamma_{0,m}$. Each
$\gamma_{0,i}$ is a semi-algebraic path parametrized by 
one of the co-ordinates $X_1,\ldots,X_k$, over some interval
$[a_{0,i},b_{0,i}]$
with $\gamma_{0,0}(a_{0,0}) = p.$
The semi-algebraic maps,
$\gamma_{0,0},\ldots,\gamma_{0,m}$ and the end-points
of their intervals of definition
$a_{0,0},b_{0,0},\ldots,a_{0,m},b_{0,m}$ are all independent of $x$
(upto the discrete choice of the path $\gamma(x)$ in $\RM(\ZZ(Q,\R^k),\{x\})$),
except $b_{0,m}$ which depends on $x_1$. 

Moreover,  $\Gamma_1(x)$ can again be decomposed into two parts, 
$\gamma_1(x)$ and $\Gamma_2(x)$
with 
$\Gamma_2(x)$ contained in 
$\ZZ(Q,\R^k)_{\bar x_2}$ and so on.

If $y = (y_1,\ldots,y_k) \in \ZZ(Q,\R^k)$ is another point such that
$x_1 \neq y_1$, then since
$\ZZ(Q,\R^k)_{x_1}$
and 
$\ZZ(Q,\R^k)_{y_1}$
are disjoint, 
it is clear that
$$
\RM(\ZZ(Q,\R^k),\{x\}) \cap \RM(\ZZ(Q,\R^k),\{y\}) = \RM(\ZZ(Q,\R^k)).
$$
Now consider a
connecting path  $\gamma(y)$ extracted from $\RM(\ZZ(Q,\R^k),\{y\})$. 
The images of $\Gamma_1(x)$ and $\Gamma_1(y)$ are
disjoint.
If the image of 
$\gamma_0(y)$ (which is contained in $\RM(\ZZ(Q,\R^k)$) follows the same
sequence of curve segments as $\gamma_0(x)$   
starting at $p$
(that is,  it consists of the same curves segments 
$\gamma_{0,0},\ldots,\gamma_{0,m}$ as in $\gamma_0(x)$),
then it is clear that the images of the 
paths $\gamma(x)$ and $\gamma(y)$ has the
property that they 
are identical upto a point and they are disjoint after it.

\hide{
(see Figure \ref{fig:figure1}).
We call this the {\em divergence property}.

\begin{center}
\begin{figure}
\label{fig:figure1}
\includegraphics[height=7cm]{figure1.pdf}
\caption{The connecting path $\Gamma(x)$}
\end{figure}
\end{center}
}
More generally, if the points $x$ and $y$ are such that,
$x_i = y_i, 1 \leq i \leq j$ and $x_{j+1} \neq y_{j+1}$, then the 
 paths $\Gamma_{j+1}(x)$ and $\Gamma_{j+1}(y)$, contained in
$\ZZ(Q,\R^k)_{\bar x_{j+1}}$
and  
$\ZZ(Q,\R^k)_{\bar y_{j+1}}$
respectively,
will be disjoint.
Moreover if the paths 
$\gamma_0(x),\ldots,\gamma_j(x)$ and 
$\gamma_0(y),\ldots,\gamma_j(y)$ 
are composed of the same
sequence of curve segments, then $\gamma(x)$ and $\gamma(y)$ will
also have the divergence property.

We now consider connecting paths in the semi-algebraic setting.
We are given  a polynomial $Q \in \R[X_1,\ldots,X_k]$
such that
$\ZZ(Q,\R^k)$ is bounded and
 a finite set of
 polynomials ${\mathcal P} \subset \D[X_1,\ldots,X_k]$
in strong $\ell$-general position with respect to $Q$.
This means that
any  $\ell+1$
polynomials belonging to 
${\mathcal P}$ have no  zeros in common with
$Q$ in $\R^k$, and 
any  $\ell$
polynomials belonging to 
${\mathcal P}$ have at  most a finite number of  zeros in common with
$Q$ in $\R^k$. 

For every point $x$ of $\ZZ(Q,\R^k)$, we denote by 
 $\sigma(x)$ the sign condition
on ${\mathcal P}$ at $x$. Let
$
 {\RR}(\overline\sigma(x),\ZZ(Q,\R^k)) =
\{x\in\ZZ(Q,\R^k)\mid\bigwedge_{P \in {\mathcal P}}\;\; \sign(P(x)) \in
\overline\sigma(x)(P)\},
$
where $\overline\sigma$ is the relaxation of $\sigma$
 defined by
$$
\left\{\begin{array}{ccc} \overline\sigma =
\{0\} &\mbox{ if }&\sigma =0,\\
\overline\sigma = \{0,1\} &\mbox{ if
}&\sigma =1,\\
\overline\sigma = \{0,-1\}& \mbox{ if }&\sigma
=-1.
\end{array}\right.
$$
We say that
$\overline \sigma(x)$ is the weak sign condition defined by
 $x$ on ${\mathcal P}$. We denote by ${\mathcal P}(x)$
the union of $\{Q\}$ and the set of
polynomials in ${\mathcal P}$ vanishing at $x$.

The connecting algorithm associates to
$x \in \ZZ(Q,\R^k)$ a path entirely 
contained in the realization of 
$\overline \sigma(x)$
connecting $x$ to a distinguished point
of the roadmap of some 
 $\ZZ({\mathcal P}',\R^k)$,
with ${\mathcal P}(x) \subset {\mathcal P}' .$ 
The connecting algorithm proceeds as follows:
construct a path $\gamma$ connecting
 a distinguished point of 
 $\RM(\ZZ(Q,\R^k))$ 
to $x$
contained in   $\RM(\ZZ(Q,\R^k),\{x\})$.
If no polynomial of ${\mathcal P} \setminus {\mathcal P}(x)$
vanishes on $\gamma$, we are done.
Otherwise let  $y$ be the 
last point of $\gamma$
such that
some polynomial of ${\mathcal P} \setminus {\mathcal P}(x)$
vanishes at $y$.
Now 
keep the part of $\gamma$ connecting $y$ to
$x$ as end of the connecting path,
and iterate the construction with $y$, noting that
the realization of
$\overline \sigma(y)$ 
is contained in the realization of
$\overline \sigma(x)$,
and  ${\mathcal P} \setminus {\mathcal P}(y)$
is in $\ell-1$ strong general position with respect to
$\ZZ({\mathcal P}(y)
,\R^k)$.

As in the algebraic case, two such connecting paths which
start with the same sequence of curve segments
will have the 
divergence
property. 
This follows from the 
divergence
property in the algebraic case 
and the recursive definition of connecting paths.

Formal descriptions and complexity analysis of the algorithms
described above for computing roadmaps and  connecting paths 
of algebraic and basic semi-algebraic sets
can be found in \cite{BPR03} (Algorithm 15.12 and Algorithm 16.8).

\section{Parametrized paths}
\label{sec:covering}
We are given  a polynomial $Q \in \R[X_1,\ldots,X_k]$
such that
$\ZZ(Q,\R^k)$ is bounded and
 a finite set of
 polynomials ${\mathcal P} \subset \D[X_1,\ldots,X_k]$
in strong $k'$-general position with respect to $Q$, where $k'$ is the
dimension of $\ZZ(Q,\R^k)$.

We show how to obtain a covering of a given 
$\mathcal P$-closed semi-algebraic set contained in $\ZZ(Q,\R^k)$
by a family of semi-algebraically contractible subsets.
The construction
is based on a parametrized version of the connecting algorithm:
we compute a family of polynomials such that for each realizable
sign condition $\sigma$ on this family, the description of the connecting
paths of different points in the realization, $\RR(\sigma,\ZZ(Q,\R^k)),$ are
uniform. 
We first define parametrized paths.
A parametrized path is a
semi-algebraic set which is a union of semi-algebraic
paths having the 
divergence property.

More precisely,
\begin{definition}
\label{def:parametrizedpath}
A parametrized path 
$\gamma$ is
 a continuous semi-algebraic mapping  from 
$V \subset \R^{k+1}
\rightarrow \R^k,$
such that,
denoting by 
$U=\pi_{1\ldots k}(V)\subset \R^k$,
there exists a semi-algebraic continuous function
$\ell: U \rightarrow [0,+\infty),$ 
and there exists  a point $a$ in $\R^k$, such that
\begin{enumerate}
\item $V =\{(x,t) \mid x \in U, 0 \le t \le \ell(x)\},$
\item $\forall \; x \in U, \; \gamma(x,0)=a$,
\item $\forall \; x \in U, \; \gamma(x,\ell(x))=x$,
\item 
$$
\displaylines{
\forall \; x \in U, \forall \; y \in U, \forall \; s \in [0,\ell(x)], \forall \; t \in [0,\ell(y)]\cr
 \left(\gamma(x,s)=\gamma(y,t) \Rightarrow s=t \right),
}
$$
\item  
$$
\displaylines{
\forall \; x \in U, \forall \; y \in U, \forall \; s \in [0,\min(\ell(x),\ell(y))]\cr
\left(\gamma(x,s)=\gamma(y,s) \Rightarrow \forall \; t \le s \; \gamma(x,t)=\gamma(y,t) \right).
}
$$
\end{enumerate}
\end{definition}

\hide{
\begin{center}
\begin{figure}
\includegraphics[height=7cm]{figure2.pdf}
\caption{A parametrized path}
\end{figure}
\end{center}
}

Given a parametrized path, $\gamma: V \rightarrow \R^k$, we will
refer to $U = \pi_{1\ldots k}(V)$ as its {\em base}.
Also,  any semi-algebraic subset 
$U' \subset U$ of the base of such a parametrized
path,
defines in a natural way  the restriction of 
$\gamma$ to the base $U'$, which is
another parametrized path,
obtained by restricting $\gamma$ to the set $V' \subset V$,
defined by 
$V' = \{(x,t) \mid x \in U', 0 \leq t \leq \ell(x) \}$.

\hide{
\begin{figure}
\epsffile{example.eps}%
\caption{Illustrative example}
\label{fig:example}
\end{figure}
}

\begin{proposition}
\label{contractible}
Let $\gamma: V \rightarrow R^k$ be a parametrized path such that
$U = \pi_{1\ldots k}(V)$ is closed and bounded.
Then, the image of $\gamma$ is semi-algebraically contractible.
\end{proposition}
We thank A. Gabrielov and N. Vorobjov for pointing out an error in a
previous version of this paper, where the same proposition was stated without
any extra condition on $U$.
In fact, Proposition \ref{contractible} is not true if we do not assume
that $U$ is closed and bounded.
\hide{
 as can be seen in the following 
example 
(see Figure 1). 
In this example $U$ is the horizontal line not containing the end-point $a$. 
The paths $\gamma$ starting from
a point on $p \in U$ consists of two curve segments -- the first one is
the curve that connects $p$ to a point $q$ on the 
other line parametrized by its $Y$ coordinate, 
followed by a straight segment joining $q$ to $a$ parametrized by its
$X$ coordinate.
The image of this parametrized path is the whole surface strictly to
the left of $a$ plus the point $a$. This is not contractible and indeed
the map $\phi$ defined in the proof below is not continuous at $(a,0)$.
This kind of phenomenon cannot happen if $U$ is closed and bounded
(see proof of Proposition \ref{contractible} below).
}

\begin{proof}(of Proposition \ref{contractible})
Let $W = \Im(\gamma)$ and $M = \sup_{x \in U}\ell(x).$
We prove that the semi-algebraic mapping 
$
\phi: W\times [0,M] \rightarrow  W
$
sending 
\begin{itemize}
\item $(\gamma(x,t),s)$ to $\gamma(x,s)$ if $t \ge s$,
\item $(\gamma(x,t),s)$ to $\gamma(x,t)$ if $t < s$.
\end{itemize}
is continuous.
Note that the map $\phi$ is well-defined,
since
$\gamma(x,t) = \gamma(x',t') \Rightarrow t = t',$
by condition (4).

Since $\phi$ satisfies 
$$
\displaylines{
\phi(\gamma(x,t),0)=a, \cr
\phi(\gamma(x,t),M)=\gamma(x,t),
}
$$
this gives a semi-algebraic continuous contraction from  $W$ to $\{a\}$.

Let $w \in W, s \in [0,M]$. 
Let $\varepsilon > 0$ be an infinitesimal, and
let $(w',s') \in \E(W \times [0,M],\R\la\eps\ra)$ be such that
$\lim_{\varepsilon}(w',s') = (w,s)$.
In order to prove the continuity of $\phi$ at $w$ it suffices to  
prove that $\lim_{\varepsilon}\E(\phi,\R\la\eps\ra)(w',s') = 
\phi(w,s).$

Let $w = \gamma(x,t)$ for some $x \in U, t \in [0,\ell(x)],$ 
and similarly let $w' = (x',t')$ 
for some $x' \in \E(U,\R\la\eps\ra)$ and 
$t' \in [0,\E(\ell,\R\la\eps\ra)(x')].$
Note that $\lim_{\varepsilon}(x') \in U$ since $U$ is closed and bounded
and $\lim_\eps t' \in [0,\ell(\lim_\eps x')].$

Now, 
$$
\begin{array}{lll}
\gamma(x,t) &=& w \\
&=& \lim_{\eps}(w') \\
&=& \lim_{\eps}\E(\gamma,\R\la\eps\ra)(x',t')\\
&=& \gamma(\lim_\eps x',\lim_\eps t').
\end{array}
$$ 
Condition (4) now implies that $\lim_{\eps} t' = t.$ 

Without loss of generality let $t' \geq t$. The other case is symmetric.  
We have the following two sub-cases.

Case $s' > t'$: 
Since $s,t \in \R$ and $\lim_\eps s' = s$
and $\lim_{\varepsilon} t' = t$, we must have that
$s \geq t.$
In this case $\E(\phi,\R\la\eps\ra)(w', s') = 
\E(\gamma,\R\la\eps\ra)(x',t').$
Then, 
$$
\begin{array}{lll}
\lim_{\varepsilon} \E(\phi,\R\la\eps\ra)(w',s') &=& 
\lim_{\varepsilon}\E(\gamma,\R\la\eps\ra)( x',t') \\
&=& 
\lim_{\varepsilon} w'  \\
&=& w \\
&=& \phi(w,s).
\end{array}
$$

Case $s' \leq  t'$: 
Again, since $s,t \in \R$ and  $\lim_{\varepsilon} s' = s$ and
$\lim_{\varepsilon} t' = t$,
we must have that $s \leq t$.

In this case we have,
$$
\begin{array} {lll}
\lim_{\varepsilon} \phi(w',s') &=& 
\lim_{\varepsilon}\E(\gamma,\R\la\eps\ra)( x',s') \\
&=& 
\gamma(\lim_{\varepsilon} x',\lim_{\varepsilon} s') \\
&=& \gamma(\lim_{\varepsilon} x',s).
\end{array}
$$

Now, 
$$
\begin{array}{lll}
\gamma(\lim_{\varepsilon} x', t) &=&
\gamma(\lim_{\varepsilon} x', \lim_{\varepsilon} t') \\
&=& \lim_{\varepsilon} \E(\gamma,\R\la\eps\ra)(x',t') \\
&=& \lim_{\varepsilon} w' \\
&=& w \\
&=& \gamma(x,t).
\end{array}
$$

Thus, by condition (5) we have that $\gamma(\lim_{\varepsilon} x',s'')
= \gamma(x,s'')$ for all $s'' \leq t$.
Since, $s \leq t$, this implies,
$$
\begin{array}{lll}
\lim_{\varepsilon} \E(\phi,\R\la\eps\ra)(w',s') &=& 
\lim_\eps \E(\gamma,\R\la\eps\ra)(w',s') \\
&=&\gamma(\lim_{\varepsilon} x',\lim_\eps s') \\
&=& \gamma(x,s)\\
&=& \phi(w,s).
\end{array}
$$
This proves the continuity of $\phi$.
\end{proof}

We now describe how to compute parametrized paths 
in single exponential time
using a parametrized version of the connecting algorithm.
We describe the input, output and
complexity of the algorithm which appears in \cite{BPR03} (Algorithm 16.15).

\begin{algorithm}(Parametrized Bounded Connecting)
\label{16:alg:paramconnecting}
\begin{itemize}
\item[] {\sc Input.}
\begin{itemize}
\item a polynomial $Q \in \D[X_1,\ldots,X_k]$,
such that $\ZZ(Q,\R^k)\subset B(0,1/c)$,
\item  a finite set of
 polynomials ${\mathcal P}
\subset \D[X_1,\ldots,X_k]$
in strong $k'$-general position with respect to
$Q$, where $k'$ is the dimension of $\ZZ(Q,\R^k)$.
\end{itemize}
\item[] {\sc Output.}
\begin{itemize}
\item
a finite set of polynomials ${\mathcal A} \subset \R[X_1,\ldots,X_k]$,
\item 
a finite set $\Theta$ of quantifier free formulas, 
with atoms of the form $P = 0, P > 0, P< 0, \;P \in {\mathcal A}$,
such that for every semi-algebraically connected component $S$ of
the realization of
every weak sign condition on ${\mathcal P}$ on $\ZZ(Q, \R^k)$,
there exists a subset $\Theta(S) \subset \Theta$
 such that
$
\displaystyle{
S=\bigcup_{\theta \in \Theta(S)} \RR(\theta,\ZZ(Q,\R^k)),
}
$
\item for every $\theta \in \Theta$,
a parametrized path
$$
\displaylines{
\gamma_\theta : V_\theta \rightarrow \R^k,
}
$$
with base $U_\theta =  \RR(\theta,\ZZ(Q,\R^k))$, 
such that for each $y \in \RR(\theta,\ZZ(Q,\R^k)),$ 
$\Im \;\gamma_\theta(y,\cdot)$  is a semi-algebraic path which
connects the point $y$
to a distinguished point  $a_\theta$
of some roadmap
$\RM(\ZZ({\mathcal P}' \cup \{Q\}, \R^k))$
where ${\mathcal P}'\subset {\mathcal P}$, staying inside
${\RR}(\overline\sigma(y),\ZZ(Q,\R^k)).$

\end{itemize}

\item[] {\sc Complexity.}
$s^{k'+1} d^{O(k^4)}$, where
 $s$ is a bound on the number of elements of
${\mathcal  P}$
and $d$ is a bound on the degrees of $Q$ and
the elements of ${\mathcal  P}$.
\end{itemize}
\end{algorithm}

\vspace{.1in}
\noindent
{\sc Proof of correctness.}
Given the proof of correctness of Algorithm 16.15 (Parametrized Bounded Connecting)
in \cite{BPR03},
the 
only extra property 
that we need to prove is that
for each $\theta \in \Theta$ ,
$\gamma_\theta$ is a parametrized path. 
It is easy to see that 
$\gamma_\theta$
satisfies the conditions of Definition \ref{def:parametrizedpath},
using the
divergence property of the paths $\gamma(y,\cdot)$
(see discussion in  Section \ref{15:sec:desc}).
\eop

\section{Constructing coverings of closed semi-algebraic sets
by closed contractible sets}
\label{sec:acycliccov}
We are again given  a polynomial $Q \in \R[X_1,\ldots,X_k]$
such that
$\ZZ(Q,\R^k)$ is bounded and
 a finite set of
 polynomials ${\mathcal P} \subset \D[X_1,\ldots,X_k]$
in strong $k'$-general position with respect to $Q$, where $k'$ is the
dimension of $\ZZ(Q,\R^k)$.
We describe 
an algorithm for computing closed contractible coverings of
${\mathcal P}$-closed semi-algebraic sets,
using the results of Section \ref{sec:covering}.

For the rest of this section we fix a
${\mathcal P}$-closed semi-algebraic set $S$ contained in
$\ZZ(Q,\R^k)$ and let $\#{\mathcal A} = t.$
We denote by ${\rm Sign}({\mathcal A},S)$
the set of realizable sign conditions 
of ${\mathcal A}$ on $\ZZ(Q,\R^k)$ whose realizations are  contained in $S$.
We continue to 
follow the notations of Algorithm \ref{16:alg:paramconnecting}.
For each $\sigma \in {\rm Sign}({\mathcal A},S)$ 
$\RR(\sigma,\ZZ(Q,\R^k))$ is contained in 
$\RR(\theta,\ZZ(Q,\R^k))$ for some $\theta \in \Theta$.
We denote by $\gamma_\sigma$ 
the restriction of $\gamma_\theta$ to the base $\RR(\sigma,\ZZ(Q,\R^k))$.
Since $\RR(\sigma,\ZZ(Q,\R^k))$ is not necessarily closed and bounded,
$\Im\; \gamma_\sigma$ might not be contractible. In order to ensure
contractibility, we restrict the base of  $\gamma_\sigma$
to a slightly smaller set which is closed, using infinitesimals.

We  introduce infinitesimals
$\varepsilon_{2t} \gg\varepsilon_{2t-1} \gg \cdots \gg \varepsilon_2 \gg 
\varepsilon_1 >0$. 
For $i=1,\ldots,2t$ we will denote by $\R_i$ the field
$\R\langle\varepsilon_{2t}\rangle\cdots\langle\varepsilon_{i}\rangle$ 
and denote by
$\R'$ the field $\R_1$. 

For $\sigma \in {\rm Sign}({\mathcal A},S)$ we define the
level of $\sigma$ by,
\[
{\rm level}(\sigma)= 
\#\{P \in {\mathcal A}\;\mid\; \sigma(P) = 0 \}.
\]

Given $\sigma \in {\rm Sign}(\mathcal A,S)$, with 
${\rm level}(\sigma) = j,$ we denote by
$\RR(\sigma_-)$ the set defined on $\ZZ(Q,\R_{2j}^k)$ by
the formula $\sigma_-$ obtained by taking the 
conjunction of
$$
\begin{array}{l}
P = 0,  \mbox{ for each } P \in {\mathcal A}
\mbox{ such that } \sigma(P) = 0, \cr
P \geq \eps_{2j},  \mbox{ for each } P \in {\mathcal A}
\mbox{ such that } \sigma(P) = 1, \cr
P \leq -\eps_{2j}, \mbox{ for each }  P \in {\mathcal A} 
\mbox{ such that } \sigma(P) = -1.
\end{array}
$$
Notice that $\RR(\sigma_-)$ is closed, bounded and contained in 
$\RR(\sigma,\ZZ(Q,\R_{2j}^k))$. 
\hide{
We will denote by $\gamma_{\sigma_-}$
the restriction of $\gamma_{\sigma}$ to the base
$\RR(\sigma_-)$.
}
Proposition \ref{contractible} implies,

\begin{proposition}
\label{prop:contractible2}
$\gamma_\sigma(\RR(\sigma_-))$
 is semi-algebraically contractible.
\end{proposition}

Note that the sets 
$\gamma_\sigma(\RR(\sigma_-))$
do not necessarily
cover $S$. So we are going to enlarge them, preserving
contractibility, to obtain a covering of $S$.

Given $\sigma \in {\rm Sign}(\mathcal A,S)$, with 
${\rm level}(\sigma) = j,$ we denote by
$\RR(\sigma_-^+)$ the set defined on $\ZZ(Q,\R_{2j-1}^k),$ 
by the formula $\sigma_-^+$ obtained by taking the 
conjunction of
$$
\begin{array}{l}
-\eps_{2j-1}\le P \le \eps_{2j-1},  \mbox{ for each } P \in {\mathcal A}
\mbox{ such that } \sigma(P) = 0, \cr
P \geq \eps_{2j},  \mbox{ for each } P \in {\mathcal A}
\mbox{ such that } \sigma(P) = 1, \cr
P \leq -\eps_{2j}, \mbox{ for each }  P \in {\mathcal A} 
\mbox{ such that } \sigma(P) = -1.
\end{array}
$$
with the formula $\phi$ defining $S$.
Let 
$C_\sigma$ be the set defined by,
$$
C_\sigma=\gamma_\sigma(\RR(\sigma_-))\cup \RR(\sigma_-^+).
$$

We now prove that 
\begin{proposition}
\label{prop:contractible3}
 $C_\sigma$ is semi-algebraically contractible.
\end{proposition}

Let $C$ be a closed and bounded semi-algebraic set
contained in $\R\la \eps\ra^k$.
We can suppose without loss of generality 
that $C$ is defined over $\R[\eps]$
(see for example Proposition 2.80 in \cite{BPR03}).
We denote by $C_t$ the semi-algebraic subset of 
$\R^k$ defined by replacing $\eps$ by $t$ in the definition of
$C$. Note that $C_\eps$ is nothing but $C$.

We are going to use the following lemma.

\begin{lemma}\label{lem:contraction}
Let $B$ be a closed and bounded semi-algebraic set
contained in $\R^k$
and let $C$ be a closed and bounded semi-algebraic set
contained in $\R\la \eps\ra^k$.
If there exists $t_0$ such that for every $t<t'<t_0$, $C_t \subset C_{t'}$ 
and
$\lim_\eps(C)=B$,
then $\E(B,R\la \eps\ra)$ 
has the same homotopy type as
$C$.
\end{lemma}
\begin{proof}
Hardt's Triviality theorem implies 
that there exists $t_0 > 0$, 
and a homeomorphism
\[
\phi_{t_0} : C_{t_0} \times (0,t_0] \rightarrow \cup_{0<t\le t_0} C_{t}
\]
which preserves $C_{t_0}$.
Replacing $t_0$ by $\eps$ gives a homeomorphism
\[
\phi_{\eps} : C \times (0,\eps] \rightarrow \cup_{0<t\le \eps} C_{t}.
\]
Defining
\[
\psi : C \times [0,\eps] \rightarrow C
\]
by 
$$\begin{cases}
\psi(x,s) = \pi_{1\ldots k}\circ \phi(x,s),\mbox{ if } s>0\cr
\psi(x,0) = \lim_{s \rightarrow 0_+}\pi_{1\ldots k}\circ \phi(x,s),       
\end{cases}$$
it is clear that $\psi$ is 
a semi-algebraic retraction of $C$
to $\E(B,\R\la \eps \ra).$ 
\end{proof}

We now prove Proposition \ref{prop:contractible3}.
\begin{proof}(of Proposition \ref{prop:contractible3})
Apply Lemma \ref {lem:contraction} to $C_\sigma$
and 
$\E(\gamma_\sigma(\RR(\sigma_-)),\R_{2j-1})$: 
thus 
$C_\sigma$ can be semi-algebraically retracted to
$\E(\gamma_\sigma(\RR(\sigma_-)),\R_{2j-1})$. 

Since 
$\E(\gamma_\sigma(\RR(\sigma_-)),\R_{2j-1})$ 
is semi-algebraically
contractible, so is $C_\sigma$.
\hide{
Replacing $\eps_{2j-1}$ in the definition of  $C_\sigma$ by
a new variable $T$, we obtain a closed and bounded semi-algebraic set 
$B_\sigma \subset \R_{2j}^{k+1}$. For $t \in \R_{2j}$, we will denote
by $B_\sigma(t) = \pi_{k+1}^{-1}(t) \cap B_\sigma.$
Hardt's Triviality theorem implies 
that there exists $t_0 > 0$, such that for all $t \in [0,t_0]$
there exists a semi-algebraic retraction of $\pi_{1\ldots k}(B_\sigma(t))$
to $\pi_{1\ldots k}(B_\sigma(0)) = \Im (\gamma_{\sigma_-}).$
Now, replacing $t_0$ by $\epsilon_{2j-1}$ we have that
$C_\sigma$ can be semi-algebraically retracted to
$\E(\Im (\gamma_{\sigma_-}),\R_{2j-1})$. 
But since $\E(\Im(\gamma_{\sigma_-}),\R_{2j-1})$ 
is semi-algebraically
contractible, so is $C_\sigma$.
}
\end{proof}

We now prove that the sets $\E(C_\sigma,\R')$ form a covering of 
$\E(S,\R')$.
\begin{proposition}(Covering property)
\label{prop:covering}
\[
\E(S,\R') = \bigcup_{\sigma \in {\rm Sign}({\mathcal A},S)} \E(C_\sigma,\R').
\]
\end{proposition}

The proposition is an immediate consequence of the following stronger result.
\begin{proposition}
\label{prop:coveringstrong}
\[
\E(S,\R') = \bigcup_{\sigma \in {\rm Sign}({\mathcal A},S)} \RR(\sigma_-^+,\R'^k).
\]
\end{proposition}
\begin{proof}
By definition,
$
\displaystyle{
\E(S,\R') \supset 
\bigcup_{\sigma \in {\rm Sign}({\mathcal A},S)} \RR(\sigma_-^+,\R'^k).
}
$
We now prove the reverse inclusion.
Clearly, we have that
$
\displaystyle{
S = \bigcup_{\sigma \in {\rm Sign}({\mathcal A},S)} \RR(\sigma,\R).
}
$ 
Let $x \in \E(S,\R')$ and $\sigma$ be the sign condition of the 
family ${\mathcal A}$ at $x$ and let
${\rm level}(\sigma) = j$. If 
$x \in \RR(\sigma_-^+,\R'^k)$, we are done.
Otherwise, there exists $P \in \mathcal A,$ such that
$x$ satisfies either 
$0< P(x)<\eps_{2j}$ 
or
$-\eps_{2j} < P(x)<0.$
Let 
${\mathcal B} = \{P \in {\mathcal A}\;\mid\; \lim_{\eps_{2j}} P(x) = 0\}.$ 
Clearly $\#{\mathcal B}=j' > j$.
Let $y = \lim_{\eps_{2j}} x$. Since, $\E(S,\R')$ is 
closed and bounded and $x \in \E(S,\R')$,
$y$ is also in $\E(S,\R')$. 
Let $\tau$ be the sign condition of ${\mathcal A}$
at $y$ with ${\rm level}(\tau) = j' > j$. 
If $x \in \RR(\tau_-^+,\R'^k)$ we are done.
Otherwise, for every $P \in \mathcal A$ such that $P(y)=0$, 
we have that 
$-\eps_{2j'-1}\le P(x)\le \eps_{2j'-1}$,
since $\lim_{\eps_{2j}}(P(x))=P(y)=0$
and $\eps_{2j'-1}\gg \eps_{2j}.$
So
there exists $P \in \mathcal A$ such that
$x$ satisfies either 
$0 < P(x)<\eps_{2j'}$ 
or
$-\eps_{2j'} < P(x)<0,$ 
and we replace 
$\mathcal B$ by 
$ \{P \in {\mathcal A}\;\mid\; \lim_{\eps_{2j'}} P(x) = 0\},$ 
and
$y$ by $y = \lim_{\eps_{2j'}} x$.
This process must terminate after at most $t$ steps.
\end{proof}

\begin{algorithm} (Covering by Contractible Sets)
\label{16:alg:acycliccovering}
\index{Connecting!Algorithm!Acyclic Covering}
\begin{itemize}
\item[]
\item []{\sc Input.}
\begin{itemize}
\item  
a finite set of $s$ polynomials ${\mathcal P} \subset \D[X_1,\ldots,X_k]$
in strong $k$-general position on $\R^k,$
with $\deg(P_i) \leq d$ for $1 \leq i \leq s,$
\item a  ${\mathcal P}$-closed semi-algebraic set $S$, 
 contained  in the sphere of center $0$ and radius $r$,
defined by a ${\mathcal P}$-closed formula $\phi$.
\end{itemize}
\item[] {\sc Output.}
a set of formulas $\{\phi_1,\ldots,\phi_M\}$ such that
\begin{itemize}
\item 
each  $\RR(\phi_i,\R'^k)$ is 
semi-algebraically contractible, and
\item 
$\displaystyle{
\bigcup_{1 \leq i \leq M} 
\RR(\phi_i,\R'^k) = \E(S,\R').
}
$
\end{itemize}

\item []{\sc Procedure.}
\begin{itemize}
\item[Step 1]
Let $Q=X_1^2+\ldots + X_k^2-r^2$.
Call Algorithm \ref{16:alg:paramconnecting}
(Parametrized Bounded Connecting)
with input $Q, {\mathcal P}$.
Let ${\mathcal A}$ be the family of polynomials output.
\item[Step 2]
Compute the set of realizable sign conditions 
${\rm Sign}({\mathcal A},S)$ 
using
Algorithm 13.37 (Sampling on an Algebraic Set) in \cite{BPR03}.
\item[Step 3]
Using Algorithm 14.21 (Quantifier Elimination) in \cite{BPR03},
eliminate one variable to compute the image of the semi-algebraic map 
$\gamma_{\sigma_-}$. 
Finally, output the set of formulas 
$
\{ \phi_\sigma \;\mid\; \sigma \in {\rm Sign}({\mathcal A},S)\}
$
describing 
the semi-algebraic set $C_\sigma$.
\end{itemize}
\item[] {\sc Complexity.}
The complexity of the algorithm is bounded by $s^{(k+1)^2}d^{O(k^5)}$.
\end{itemize}
\end{algorithm}

\vspace{.1in}
\noindent
{\sc Proof of correctness.}
The correctness of the algorithm is a consequence of 
Propositions \ref{prop:contractible3}
and \ref{prop:covering} and the correctness of
Algorithm \ref{16:alg:paramconnecting}
(Parametrized Bounded Connecting), as well as 
the correctness of Algorithms 13.37 and 14.20  in \cite{BPR03}.
\eop

\vspace{.1in}
\noindent
{\sc Complexity analysis.}
The complexity of Step 1 of the algorithm is bounded by
$s^{k+1} d^{O(k^4)}$, where
$s$ is a bound on the number of elements of
${\mathcal  P}$
and $d$ is a bound on the degrees of 
the elements of ${\mathcal  P}$,
using the complexity analysis of Algorithm 
\ref{16:alg:paramconnecting} (Parametrized Bounded Connecting).
The number of polynomials in ${\mathcal A}$ is $s^{k+1} d^{O(k^4)}$
and their degrees are bounded by  $d^{O(k^3)}$. Thus 
the complexity of 
computing ${\rm Sign}({\mathcal A},S)$
is bounded by
$s^{(k+1)^2}d^{O(k^5)}$ using 
Algorithm 13.37 (Sampling on an Algebraic Set) in \cite{BPR03}.
In Step 3 of the algorithm there is a call to Algorithm 14.21 (Quantifier 
Elimination) in \cite{BPR03}. There are two blocks of variables of size
$k$ and $2$ respectively. The number and degrees of the input polynomials 
are bounded by  $s^{k+1} d^{O(k^4)}$ and 
$d^{O(k^3)}$ respectively. Moreover, observe that
even though we introduced $2s$ infinitesimals,
each arithmetic operation is
performed in the ring $\D$ adjoined with at most $O(k)$ infinitesimals
since the polynomials $\{P,P\pm \eps_{2j},P\pm \eps_{2j-1}, P \in {\mathcal P},
1 \leq j \leq s\}$ are in strong general position.
Thus, the complexity of this step is bounded by
$s^{(k+1)^2}d^{O(k^5)}$ using the complexity analysis of  
Algorithm 14.21 (Quantifier Elimination) in \cite{BPR03} and the fact that
each arithmetic operation costs at most $d^{O(k^5)}$ in terms of
arithmetic operations in the ring $\D$. 
\eop

\section{Topological Preliminaries}
\label{sec:top}
We first recall some results from algebraic topology which enable us to
compute the first Betti number of any given  closed and bounded semi-algebraic
set, from the inclusion relationships amongst the connected components
of the various pair-wise and triple-wise intersections of  the elements
of a covering of the given set by a finite number of contractible sets.

\subsection{Generalized Mayer-Vietoris exact sequence}
Let $A_1,\ldots,A_n$ be sub-complexes of a finite simplicial complex
$A$ such that $A = A_1 \cup \cdots \cup A_n$. 
Note that the intersections of any number of the sub-complexes, $A_i$,
is again a sub-complex of $A$.
We will denote by $A_{i_0,\ldots,i_p}$ the sub-complex
$A_{i_0} \cap \cdots \cap A_{i_p}$.

Let $C^i(A)$ denote the
$\Q$-vector space of $i$ co-chains of $A$, and 
$C^{\bullet}(A),$ the complex
\[ \cdots \rightarrow C^{q-1}(A) \stackrel{d}{\longrightarrow}
C^q(A) \stackrel{d}{\longrightarrow} C^{q+1}(A) \rightarrow \cdots
\]
where 
$d: C^q(A) \rightarrow C^{q+1}(A)$ 
are the usual co-boundary homomorphisms. 
More precisely, given $\omega \in  C^q(A)$, and a 
$q+1$ simplex $[a_0,\ldots,a_{q+1}] \in A$,
\begin{equation}
\label{co-boundary}
d\omega([a_0,\ldots,a_{q+1}]) = \sum_{0 \leq i \leq q+1}
            (-1)^i \omega([a_0,\ldots, \hat{a_i}, \ldots,a_{q+1}])
\end{equation}
(here and everywhere else in the paper $\hat{}$ denotes
omission). Now extend $d\omega$ to a linear form on 
all of $C_{q+1}(A)$ by linearity, to
obtain an element of $C^{q+1}(A).$ 

The generalized Mayer-Vietoris sequence is the following:
\[
0 \longrightarrow C^{\bullet}(A) \stackrel{r}{\longrightarrow}
\prod_{i_0} C^{\bullet}(A_{i_0}) 
\stackrel{\delta_1}{\longrightarrow}
\prod_{i_0 < i_1} C^{\bullet}(A_{i_0,i_1}) 
\]
\[
\cdots \stackrel{\delta_{p-1}}{\longrightarrow}
\prod_{i_0 < \cdots < i_p}C^{\bullet}(A_{i_0,\ldots,i_p})
 \stackrel{\delta_p}{\longrightarrow}
\prod_{i_0 < \cdots < i_{p+1}}C^{\bullet}(A_{i_0,\ldots,i_{p+1}})
\cdots
\]
where $r$ is induced by restriction and the connecting homomorphisms 
$\delta$ are described below.

Given an $\omega \in \prod_{i_0 < \cdots < i_p}C^q(A_{i_0,\ldots,i_p})$ 
we define $\delta(\omega)$ as follows:\\
first note that 
$\delta(\omega) \in \prod_{i_0 < \cdots < i_{p+1}}C^q(A_{i_0,\ldots,i_{p+1}})$, and it suffices to define
$\delta(\omega)_{i_0,\ldots,i_{p+1}}$ for each $(p+2)$-tuple
$0 \leq i_0 < \cdots < i_{p+1} \leq n$.
Note that, $\delta(\omega)_{i_0,\ldots,i_{p+1}}$ is a linear
form on the vector space, $C_q(A_{i_0,\ldots,i_{p+1}})$, and
hence is determined by its values on the $q$-simplices 
in the complex $A_{i_0,\ldots,i_{p+1}}$. Furthermore, 
each $q$-simplex, $s \in  A_{i_0,\ldots,i_{p+1}}$ is automatically
a simplex of the complexes 
$A_{i_0,\ldots,\hat{i_i},\ldots i_{p+1}}, \; 0 \leq i \leq p+1.$

We define,
\[
\label{delta}
(\delta \omega)_{i_0,\ldots,i_{p+1}}(s) =  
 \sum_{0 \leq i \leq p+1} (-1)^i \omega_{i_0,\ldots,\hat{i_i},\ldots,i_{p+1}}(s),
\]
 (here and everywhere else in the paper $\hat{}$ denotes
omission).
The fact that the generalized Mayer-Vietoris  sequence
is exact is classical (see \cite{B03} for example).

The cohomology groups $H^0(A_{i_0,\ldots,i_p})$ 
are isomorphic to the $\Q$-vector space of locally
constant functions on $A_{i_0,\ldots,i_p}$ and the induced homomorphisms,
\[
\delta_p: 
H^*(A_{i_0,\ldots,i_p}) \rightarrow  H^*(A_{i_0,\ldots,i_{p+1}})
\]
are then given by generalized restrictions, i.e.
for 
$$
\displaylines{
\phi \in \oplus_{1 \leq i_0 < \cdots < i_p \leq n}H^0(A_{i_0,\ldots,i_p}),
}
$$ 
a locally constant function on $A_{i_0,\ldots,i_p}$, 
$$
\displaylines{
\delta_p(\phi)_{i_0,\ldots,i_{p+1}} = 
\sum_{i=0}^{p} (-1)^i 
\phi_{i_0,\ldots,\hat{i_i},\ldots,i_{p+1}}|_
{A_{i_0,\ldots,i_{p+1}}}.
}
$$

The following proposition
provides the key tool for computing the first Betti number.
\begin{proposition}
\label{prop:computebettione}
Let $A_1,\ldots,A_n$ be sub-complexes of a finite simplicial complex
$A$ such that $A = A_1 \cup \cdots \cup A_n$ and each 
$A_i$ is acyclic, that is $H^0(A_i) = \Q$ and $H^q(A_i) = 0$ for all $q >0$. 
Then, 
$b_1(A)=\dim(\Ker(\delta_2))-\dim(\Ima(\delta_1))$,
with
$$
\prod_{i}H^0(A_{i}) 
\stackrel{\delta_1}{\longrightarrow}  \prod_{i<j}
H^0(A_{i,j}) 
\stackrel{\delta_2}{\longrightarrow}  \prod_{i<j <\ell}
H^0(A_{i,j,\ell}) 
$$
\end{proposition}

To prove Proposition
\ref{prop:computebettione},
we consider the following bi-graded double complex ${\mathcal M}^{p,q}$,
with  total differential $D = \delta + (-1)^p d$, where 
\[
{\mathcal M}^{p,q} =  
\prod_{i_0,\ldots,i_p}C^q(A_{i_0,\ldots,i_p}).
\]

{\small
\[
\begin{array}{cccccccc}

& & \vdots  && \vdots  && \vdots  & \\
& &
\Big\uparrow\vcenter{\rlap{$d$}} & &
\Big\uparrow\vcenter{\rlap{$d$}} & &
\Big\uparrow\vcenter{\rlap{$d$}} &  \\

0 & \longrightarrow & \prod_{i_0}C^3(A_{i_0}) &
\stackrel{\delta}{\longrightarrow} & \prod_{i_0<i_1}
C^3(A_{i_0,i_1}) &
\stackrel{\delta}{\longrightarrow} & \prod_{i_0<i_1 <i_2}
C^3(A_{i_0,i_1,i_2}) &
 \longrightarrow 
\\

& &
\Big\uparrow\vcenter{\rlap{$d$}} & &
\Big\uparrow\vcenter{\rlap{$d$}} & &
\Big\uparrow\vcenter{\rlap{$d$}} &  \\

0 & \longrightarrow & \prod_{i_0}C^2(A_{i_0}) &
\stackrel{\delta}{\longrightarrow} & \prod_{i_0<i_1}
C^2(A_{i_0,i_1}) &
\stackrel{\delta}{\longrightarrow} & \prod_{i_0<i_1 <i_2}
C^2(A_{i_0,i_1,i_2}) &
 \longrightarrow 
\\

& &
\Big\uparrow\vcenter{\rlap{$d$}} & &
\Big\uparrow\vcenter{\rlap{$d$}} & &
\Big\uparrow\vcenter{\rlap{$d$}} &  \\

0 & \longrightarrow & \prod_{i_0}C^1(A_{i_0}) &
\stackrel{\delta}{\longrightarrow} & \prod_{i_0<i_1}
C^1(A_{i_0,i_1}) &
\stackrel{\delta}{\longrightarrow} & \prod_{i_0<i_1 <i_2}
C^1(A_{i_0,i_1,i_2}) &
 \longrightarrow 
\\

& &
\Big\uparrow\vcenter{\rlap{$d$}} & &
\Big\uparrow\vcenter{\rlap{$d$}} & &
\Big\uparrow\vcenter{\rlap{$d$}} &  \\

0 & \longrightarrow & \prod_{i_0}C^0(A_{i_0}) &
\stackrel{\delta}{\longrightarrow} & \prod_{i_0<i_1}
C^0(A_{i_0,i_1}) &
\stackrel{\delta}{\longrightarrow} & \prod_{i_0<i_1 <i_2}
C^0(A_{i_0,i_1,i_2}) &
 \longrightarrow 
\\

& &
\Big\uparrow\vcenter{\rlap{$d$}} & &
\Big\uparrow\vcenter{\rlap{$d$}} & &
\Big\uparrow\vcenter{\rlap{$d$}} &  \\

& & 0 && 0 && 0 & \\
\end{array}
\]
}

There are two spectral sequences (corresponding to taking horizontal or
vertical filtrations respectively) associated with ${\mathcal M}^{p,q}$ both
converging to $H^*_D({\mathcal M})$. The first terms of these are
${'E}_1 = H_{\delta}{\mathcal M}, {'E}_2 = H_dH_{\delta} {\mathcal M}$, and
${''E}_1 = H_d {\mathcal M}, {''E}_2 = H_\delta H_d {\mathcal M}$. 
Because of the
exactness of the generalized Mayer-Vietoris sequence, we have that,
\[
{'E}_1 =
\begin{array}{|cccccccccccccc}
& \vdots  && \vdots&& \vdots&& \vdots  && \vdots  && \\

&\Big\uparrow\vcenter{\rlap{$d$}} & &
\Big\uparrow\vcenter{\rlap{$0$}} & &
\Big\uparrow\vcenter{\rlap{$0$}} & &
\Big\uparrow\vcenter{\rlap{$0$}} & &
\Big\uparrow\vcenter{\rlap{$0$}}& & \\

& C^3(A) &
 & 0 &
 & 0 &
 & 0 &
 & 0 & \cdots &
  
\\

&\Big\uparrow\vcenter{\rlap{$d$}} & &
\Big\uparrow\vcenter{\rlap{$0$}} & &
\Big\uparrow\vcenter{\rlap{$0$}} & &
\Big\uparrow\vcenter{\rlap{$0$}} & &
\Big\uparrow\vcenter{\rlap{$0$}} & & \\

& C^2(A) &
 & 0 &
 & 0 &
 & 0 &
 & 0 & \cdots &
 
\\

&\Big\uparrow\vcenter{\rlap{$d$}} & &
\Big\uparrow\vcenter{\rlap{$0$}} & &
\Big\uparrow\vcenter{\rlap{$0$}} & &
\Big\uparrow\vcenter{\rlap{$0$}} & &
\Big\uparrow\vcenter{\rlap{$0$}} & & \\

& C^1(A) &
 & 0 &
 & 0 &
 & 0 &
 & 0 & \cdots &
 
\\

&\Big\uparrow\vcenter{\rlap{$d$}} & &
\Big\uparrow\vcenter{\rlap{$0$}} & &
\Big\uparrow\vcenter{\rlap{$0$}} & &
\Big\uparrow\vcenter{\rlap{$0$}} & &
\Big\uparrow\vcenter{\rlap{$0$}} & & \\

&C^0(A) &
 & 0 &
 & 0 &
 & 0 &
 & 0 & \cdots &
\\ \\
\hline
\end{array}
\]

and
\[
{'E}_2 = 
\begin{array}{|cccccccc}
 & \vdots  &\vdots & \vdots  & \vdots & \vdots  & \vdots\\
& &
 & &
 & &
 &  \\

 & H^3(A) & 0
 & 0 & 0
 & 0 & 0& \cdots 
\\

 &
 & &
 & &
 & & \\

  & H^2(A) & 0
 & 0 & 0
 & 0 & 0 &\cdots
 
\\

 &
 & &
 & &
 & & \\

  & H^1(A) & 0
 & 0 & 0
 & 0 & 0 &\cdots 
\\

 &
 & &
 & &
 & & \\

 & H^0(A) & 0
 & 0 & 0
 & 0 & 0 &\cdots 
\\

& 
 & &
 & &
 & & \\
\hline
\end{array}
\]

The degeneration of this sequence at $E_2$ shows that
$H^*_D({\mathcal M}) \cong H^*(A)$.

The initial term ${''E}_1$ of the second spectral sequence is given by,
{\small
\[
{''E}_1 =
\begin{array}{|ccccccc}

 & \vdots  && \vdots  && \vdots  & \\
 &
 & &
 & &
 &  \\

 & \prod_{i_0}H^3(A_{i_0}) &
\stackrel{\delta}{\longrightarrow} & \prod_{i_0<i_1}
H^3(A_{i_0,i_1}) &
\stackrel{\delta}{\longrightarrow} & \prod_{i_0<i_1 <i_2}
H^3(A_{i_0,i_1,i_2}) &
 \longrightarrow 
\\

 &
 & &
 & &
 &  \\

 & \prod_{i_0}H^2(A_{i_0}) &
\stackrel{\delta}{\longrightarrow} & \prod_{i_0<i_1}
H^2(A_{i_0,i_1}) &
\stackrel{\delta}{\longrightarrow} & \prod_{i_0<i_1 <i_2}
H^2(A_{i_0,i_1,i_2}) &
 \longrightarrow 
\\

 &
 & &
 & &
 &  \\

 & \prod_{i_0}H^1(A_{i_0}) &
\stackrel{\delta}{\longrightarrow} & \prod_{i_0<i_1}
H^1(A_{i_0,i_1}) &
\stackrel{\delta}{\longrightarrow} & \prod_{i_0<i_1 <i_2}
H^1(A_{i_0,i_1,i_2}) &
 \longrightarrow 
\\

 &
 & &
 & &
 &  \\

 & \prod_{i_0}H^0(A_{i_0}) &
\stackrel{\delta}{\longrightarrow} & \prod_{i_0<i_1}
H^0(A_{i_0,i_1}) &
\stackrel{\delta}{\longrightarrow} & \prod_{i_0<i_1 <i_2}
H^0(A_{i_0,i_1,i_2}) &
 \longrightarrow \\\\
\hline 
\end{array}
\]
}

The cohomology groups $H^0(A_{i_0,\ldots,i_p})$ occurring as
summands in the
bottom row of ${''E}_1$ are isomorphic to the $\Q$-vector space of locally
constant functions on $A_{i_0,\ldots,i_p}$ and the homomorphisms,
$''d_1: {''E}_1^{p,0} \rightarrow {''E}_1^{p+1,0}$ are 
then given by generalized restrictions, i.e.
for 
$$
\displaylines{
\phi \in \oplus_{1 \leq i_0 < \cdots < i_p \leq n}H^0(A_{i_0,\ldots,i_p}),
}
$$
with each $\phi_{i_0,\ldots,i_{p+1}}$
a locally constant function on $A_{i_0,\ldots,i_p}$, 
$$
\displaylines{
''d_1(\phi)_{i_0,\ldots,i_{p+1}} = 
\sum_{i=0}^{p} (-1)^i 
\phi_{i_0,\ldots,\hat{i_i},\ldots,i_{p+1}}|_
{A_{i_0,\ldots,i_{p+1}}}.
}
$$

\begin{proof}(Proof of Proposition \ref{prop:computebettione})
Since,  $H^q(A_{i}) = 0$ for all $q > 0$, all the terms in the first 
column of ${''E}_1$ are zero except the bottom term, and clearly
${''d}_2^{0,1} = 0$. 
Thus, ${''E}_{\infty}^{1,0} = {''E}_{2}^{1,0}$ and
${''E}_{\infty}^{0,1} = 0$. 
Thus, 
$H^1(A) \cong {''E}_{\infty}^{1,0} \oplus {''E}_{\infty}^{0,1}
 \cong {''E}_{2}^{1,0}.$
\end{proof}

\section{Computing the first Betti number in the ${\mathcal P}$-closed case}
\label{sec:closedcase}
Let $S$ be a ${\mathcal P}$-closed semi-algebraic set.
We first replace $S$ by a  ${\mathcal P}^\star$-closed and bounded semi-algebraic set,
where the elements of ${\mathcal P}^\star$ are slight modifications
of the elements of  ${\mathcal P}$,
and the family ${\mathcal P}^\star$ is in general position
and $b_i(S^\star) = b_i(S), 0 \leq i \leq k$.

Define
$$
H_i=1 + \sum_{1\leq j\leq k} i^j X_j^{d'}.
$$
where $d'$ is the smallest number strictly bigger than the degree
of all the polynomials in $\mathcal P$.
Using arguments similar to the proof
of Proposition 13.7 in 
\cite{BPR03}, it is easy to see
that the family $\mathcal P^\star$ of polynomials
 $P_i -\delta H_i, P_i +\delta H_i$, with $P_i \in \mathcal P$.
 is in 
general position in $\R \la \delta \ra^k$.

\begin{lemma}
\label{lem:star}
Denote by
$S^\star$ the set obtained by replacing any
 $P_i \ge 0$ in the definition of $S$ by  $P_i \ge -\delta H_i$ 
 and every  $P_i \le 0$ in the definition of $S$ by  $P_i \le \delta H_i$.
If $S$ is bounded, the set $\E(S,\R\la \delta \ra ^k$
is semi-algebraically homotopy equivalent to $S^\star$.
\end{lemma}
\begin{proof}
Note that $S$ is closed and bounded, $\lim_\delta S^\star=S$,and $S_t \subset S_{t'}$.
The claim follows by Lemma 
\ref{lem:contraction}.
\end{proof}

\begin{algorithm}
(First Betti Number of a ${\mathcal P}$-closed Semi-algebraic Set)
\label{alg:firstbetticlosed}
\begin{itemize}
\item[]
\item[] {\sc Input.} 
\begin{itemize}
\item  
a finite set of  polynomials ${\mathcal P} \subset \D[X_1,\ldots,X_k],$
\item a formula defining a ${\mathcal P}$-closed semi-algebraic set, $S$.
\end{itemize}

\item[]{\sc Output.}
$b_1(S).$
\item[] {\sc Procedure.}
\item[Step 1]
Let $\eps$ be an infinitesimal. 
 Replace $S$ by the semi-algebraic set $T$ defined as the intersection
 of the cylinder 
$S\times  \R\la \eps \ra$ 
with the upper hemisphere defined by
 $\eps^2(X_1^2+\ldots+X_k^2+X_{k+1}^2)=1, X_{k+1}\ge 0.$
 \item[Step 2] Replace $T$ by $T^\star$ using the notation
 of Lemma \ref{lem:star}.
\item[Step 3]
Use Algorithm \ref{16:alg:acycliccovering} (Covering by Contractible Sets)
with input $\eps^2(X_1^2+\ldots+X_k^2+X_{k+1}^2)-4$ and ${\mathcal P}^\star$, 
to compute a covering of  
$T^\star$
by closed, bounded and  contractible
sets, $T_i$, described by formulas
$\phi_i$.
\item[Step 4]
Use Algorithm 16.27  (General Roadmap) in 
\cite{BPR03} to
compute exactly one sample point of
each connected component
of the pairwise and triplewise intersections of the $T_i$'s.
For every pair $i,j$ and every $k$ compute the incidence relation
between the connected components of $T^\star_{ijk}$ and $T^\star_{ij}$ as follows:
compute a roadmap of $T^\star_{ij}$, 
containing the sample points of the connected components  of $T^\star_{ijk}$
using Algorithm 16.27  (General Roadmap).
\item[Step 5]
Using linear algebra compute 
$b_1(T^\star)=\dim(\Ker(\delta_2))-\dim(\Im(\delta_1))$,
with 
$$
 \prod_{i}H^0(T^\star_{i}) 
\stackrel{\delta_1}{\rightarrow}  \prod_{i<j}
H^0(T^\star_{ij})
\stackrel{\delta_2}{\rightarrow}  \prod_{i<j <\ell}
H^0(T^\star_{ij\ell}) 
$$
\item[] {\sc Complexity.}
The complexity of the algorithm is bounded by 
$(sd)^{k^{O(1)}}$,
where $s = \#{\mathcal P}$ and $d = \max_{P \in {\mathcal P}}\deg(P).$
\end{itemize}
\end{algorithm}

\vspace{.1in}
\noindent
{\sc Proof of correctness.}
First note that $T$ is closed and bounded and has the same
Betti numbers as $S$, using the local conical structure at
infinity. It follows from Lemma \ref{lem:star}
that $T$ and $T^\star$ have the same Betti numbers. 
The correctness of the algorithm is a consequence of 
the correctness of
Algorithm \ref{16:alg:acycliccovering} (Covering by Contractible Sets),
Algorithm 16.27  (General Roadmap) 
 in 
\cite{BPR03}, and Proposition \ref{prop:computebettione}.
\eop

\vspace{.1in}
\noindent
{\sc Complexity analysis.}
\hide{
Each step is clearly singly exponential from the complexity analysis
of Algorithm \ref{16:alg:acycliccovering} (Covering by Contractible Sets), 
and
Algorithm 16.27  (General Roadmap) in 
\cite{BPR03} and the fact that the linear algebra in  Step 5 can also
be performed in singly exponential time.
}
The complexity of Step 3 of the algorithm is bounded by
$s^{(k+1)^2}d^{O(k^6)}$ using the complexity analysis of 
Algorithm \ref{16:alg:acycliccovering} (Covering by Contractible Sets) and
noticing that each arithmetic operation takes place a ring consisting
of $\D$ adjoined with at most $k$ infinitesimals.
Finally, the complexity of Step 4 is also bounded by 
$(sd)^{k^{O(1)}}$, using the complexity analysis of 
Algorithm 16.27  (General Roadmap) in \cite{BPR03}.
\eop

\section{Replacement by closed sets without changing homology}
\label{sec:GV}
In this section, we describe a construction due to Gabrielov and Vorobjov
\cite{GV}
for replacing any given semi-algebraic subset of a bounded 
semi-algebraic set by  a closed bounded semi-algebraic subset
and strengthen the result in \cite{GV} to prove that the 
new set has the same homotopy type as the original one.
Moreover, the polynomials defining the  bounded closed semi-algebraic subset
are closely related
(by infinitesimal perturbations)
to the polynomials defining the original subset. In particular, their degrees
do not increase, while the number of polynomials used in the definition
of the new set is at most twice 
the square of the number used in the definition of
the original set.

Let ${\mathcal C} \subset
\R[X_1,\ldots,X_k]$ be a finite set of polynomials with $t$ elements, and
let $S$ be a  bounded ${\mathcal C}$-closed set.
We denote by ${\rm Sign}({\mathcal C},S)$
the set of realizable sign conditions 
of ${\mathcal C}$ whose realizations are  contained in $S$.

Recall that,
for $\sigma \in {\rm Sign}({\mathcal C})$ we define the
level of $\sigma$ 
as $\#\{P \in {\mathcal C}| \sigma(P) = 0 \}.$ 
As before let,
$\varepsilon_{2t} \gg\varepsilon_{2t-1} \gg \cdots \gg \varepsilon_2 \gg 
\varepsilon_1 >0$  be infinitesimals, 
and we will denote by $\R_i$ the field
$\R\langle\varepsilon_{2t}\rangle\cdots\langle\varepsilon_{i}\rangle$ 
and denote by
$\R'$ the field $\R_1$. For $i > 2t$, $\R_i = \R$ and for
$i < 0, \R_i = \R'$. 

We now describe the construction due to Gabrielov and Vorobjov.
For each level $m, 0 \leq m \leq t$,
we denote by
${\rm Sign}_m({\mathcal C},S)$
 the subset 
of  
${\rm Sign}({\mathcal C},S)$ of elements of level $m$.

\hide{
For each sign condition $\sigma \in {\rm Sign}_m({\mathcal C},S)$, let
$$
\displaylines{
P_{\sigma} = \sum_{P \in {\mathcal C}, \sigma(P) = 0} P^2.
}
$$
}

Given $\sigma \in {\rm Sign}_m({\mathcal C},S)$ denote by 
$\RR(\sigma_+^c)$  the intersection of  $\E(S,\R_{2m})$
 with the closed semi-algebraic set defined by 
the conjunction of the inequalities,
$$
\begin{cases}
-\varepsilon_{2m} \leq P \leq \varepsilon_{2m} \mbox{ for each }  
P \in {\mathcal C} \mbox{ such that } \sigma(P) = 0,\cr
P \ge  0,  \mbox{ for each }  P \in {\mathcal C} 
\mbox{ such that } \sigma(P) = 1,\cr
P \le  0,  \mbox{ for each }  P \in {\mathcal C} 
\mbox{ such that } \sigma(P) = -1.
\end{cases}
$$
and
denote by,
$\RR({\sigma_+^o})$ the intersection of  $\E(S,\R_{2m-1})$ with 
the open semi-algebraic set defined by the conjunction of the inequalities,
$$
\begin{cases}
-\varepsilon_{2m-1} < P < \varepsilon_{2m-1} \mbox{ for each }  
P \in {\mathcal C} \mbox{ such that } \sigma(P) = 0,\cr
P >  0,  \mbox{ for each }  P \in {\mathcal C} 
\mbox{ such that } \sigma(P) = 1,\cr
P <  0,  \mbox{ for each }  P \in {\mathcal C} 
\mbox{ such that } \sigma(P) = -1.
\end{cases}
$$

Notice that,
$$
\displaylines{
\E(\RR(\sigma), \R_{2m}) \subset \RR({\sigma_+^c}), \cr
\E(\RR(\sigma), \R_{2m-1}) \subset \RR({\sigma_+^o}).
}
$$

Let $X \subset  S$ be a ${\mathcal C}$-semi-algebraic set
such that 
$\displaystyle{
X = \bigcup_{\sigma \in \Sigma}\RR(\sigma) 
}
$
with $\Sigma\subset {\rm Sign}({\mathcal C},S)$.
We denote
$\Sigma_m= \Sigma \cap {\rm Sign}_m({\mathcal C},S)$
and define a sequence of sets, 
$X^{m} \subset \R'^k$, $0 \leq m \leq t$ inductively.

\begin{itemize}
\item
Define
$X^{0} = \E(X,\R').$
\item 
For 
$0 \leq m \leq t$,
we define 
$$
\displaylines{
X^{m+1} = 
\left(
X^{m} \cup 
\bigcup_{\sigma \in \Sigma_m}\E(\RR({\sigma_+^c})
,\R')
\right)
\setminus 
\bigcup_{\sigma \in {\rm Sign}_m({\mathcal C},S)\setminus \Sigma_m}\E(\RR({\sigma_+^o}),\R')
}
$$
\end{itemize}

We denote by  $X'$ the set $ X^{t+1}.$

The following theorem is a slight 
strengthening  of a result in \cite{GV} (where it is shown that the
sum of the Betti numbers of $X$ and $X'$ are equal) 
and the proof is very close in spirit to the one in \cite{GV}.

\begin{theorem}
\label{the:GV}
The sets $\E(X,R')$ and $X'$ are  semi-algebraically 
homotopy equivalent. In particular,
$$
\displaylines{
H_*(X) \cong H_*(X').
}
$$
\end{theorem}

For the purpose of the proof we introduce several new 
families of sets defined inductively.

For each $p,\; 0 \leq p \leq t+1$ we define sets, 
$Y_{p} \subset \R_{2p}^k, Z_{p} \subset \R_{2p-1}^k$ as follows.

\begin{itemize}
\item
We define
$$
\displaylines{
Y_{p}^{p} =  \E(X,\R_{2p}) \cup  \bigcup_{\sigma \in \Sigma_p}
\RR({\sigma_+^c}), \cr
Z_{p}^{p} = 
\E(Y_{p}^{p},\R_{2p-1}) \setminus
\bigcup_{\sigma \in {\rm Sign}_p({\mathcal C},S)\setminus \Sigma_p}\RR({\sigma_+^o}).
}
$$

\item For 
$p \leq m \leq t$,
we define 
$$
\displaylines{
Y_{p}^{m+1} = 
\left(
Y_{p}^{m} \cup \bigcup_{\sigma \in \Sigma_m}
\E(\RR({\sigma_+^c}),\R_{2p})
\right)
\setminus
\bigcup_{\sigma \in {\rm Sign}_m({\mathcal C},S)\setminus \Sigma_m}
\E(\RR({\sigma_+^o}),\R_{2p})
 \cr
Z_{p}^{m+1} = 
\left(
Z_{p}^{m} \cup  \bigcup_{\sigma \in \Sigma_m}
\E(\RR({\sigma_+^c}),\R_{2p-1})
\right)
\setminus
\bigcup_{\sigma \in {\rm Sign}_m({\mathcal C},S)
\setminus \Sigma_m}\E(\RR({\sigma_+^o}),\R_{2p-1}).
}
$$
\end{itemize}

We denote by $Y_{p} \subset \R_{2p}^k$  
(respectively, $Z_p \subset \R_{2p-1}^k$) 
the set $Y_{p}^{t+1}$ (respectively, $Z_p^{t+1}$).

Note that
\begin{itemize}
\item
$X = Y_{t+1} = Z_{t+1},$ and 
\item 
$Z_0 = X'.$ 
\end{itemize}

Notice also that for each $p, 0 \leq p \leq t,$
\begin{enumerate}
\item
$\E(Z_{p+1}^{p+1},\R_{2p}) \subset Y_p^p  ,$
\item
$Z_{p}^p \subset \E(Y_{p}^p, \R_{2p-1})  .$
  \end{enumerate}
 
The following inclusions now follow directly from the definitions
of $Y_p$ and $Z_p$.

\begin{lemma}
\label{lem:inclusions}
For each $p, 0 \leq p \leq t,$
\begin{enumerate}
\item
$\E(Z_{p+1},\R_{2p}) \subset Y_p,$
\item
$Z_{p} \subset \E(Y_{p}, \R_{2p-1})  .$
  \end{enumerate}
\end{lemma}

\hide{
\begin{lemma}
\label{lem:obvious}
Let $A \subset B$. Then $(A \cap X)\cup Y \subset (B \cap X)\cup Y$.
\end{lemma}
}
\hide{
In each case, the sets involved are constructed by applying an identical
sequence of set operations (that is adding or subtracting the same sequence
of sets) to a pair of initial sets which satisfy an obvious 
inclusion relationship. The claim follows once we note that the inclusion
relationship is preserved under the operation of
adding or subtracting the same set to both sides.
}

We now prove that in both the inclusions in Lemma \ref{lem:inclusions}
above,
the pairs of sets are in fact semi-algebraically homotopy equivalent.
These suffice to prove Theorem  \ref{the:GV}.

\begin{lemma}
\label{lem:2}
For $1 \leq p \leq t$, $Y_p$ is semi-algebraically homotopy equivalent to 
$\E(Z_{p+1},\R_{2p}).$
\end{lemma}

\begin{proof}
Let $Y_{p}(u) \subset \R_{2p+1}^k$ denote set obtained by replacing the
infinitesimal $\varepsilon_{2p}$ in the definition of $Y_{p}$ by $u$,
and for $u_0 > 0$, we will denote by
Let $Y_p((0,u_0]) \subset \R_{2p+1}^{k+1}$  the set 
$\{(x,u)| x \in Y_{p}(u), u \in (0,u_0] \}.$ 

By Hardt's triviality theorem
there exist $u_0 \in \R_{2p+1}, \;u_0 > 0$ and a homeomorphism,
\[
\psi: Y_{p}(u_0) \times (0,u_0] \rightarrow 
Y_p((0,u_0]),
\]
such that 
\begin{enumerate}
\item
$\pi_{k+1}(\phi(x,u)) = u,$ 
\item
$\psi(x,u_0)  = (x,u_0)$ for $x \in Y_{p}(u_0),$ 
and
\item
for all $u \in (0,u_0],$
and for every sign condition $\sigma$ of the
family, $\cup_{P \in {\mathcal C}} \{P, P\pm\varepsilon_{2t},\ldots,P\pm
\varepsilon_{2p+1}\},$
$\psi(\cdot,u)$ restricts to a homeomorphism of 
$\RR(\sigma,Y_p(u_0))$ to $\RR(\sigma,Y_p(u)).$
\end{enumerate}

Now, we specialize $u_0$ to $\eps_{2p}$ and denote the map corresponding
to $\psi$ by $\phi$.

For $\sigma \in \Sigma_p,$ we define,
$\RR(\sigma_{++}^o)$ to be the set defined by,
$$
\begin{cases}
-2\eps_{2p}  < P <  2\eps_{2p},  \mbox{ for each }  P \in {\mathcal C} 
\mbox{ such that } \sigma(P) = 0,\cr
P >  -\eps_{2p},  \mbox{ for each }  P \in {\mathcal C} 
\mbox{ such that } \sigma(P) = 1,\cr
P <  \eps_{2p},  \mbox{ for each }  P \in {\mathcal C} 
\mbox{ such that } \sigma(P) = -1.
\end{cases}
$$

Let 
$ \lambda: Y_p \rightarrow \R_{2p}$ 
be a semi-algebraic continuous function such that,
$$
\begin{cases}
\lambda(x) = 1,\mbox{ on } Y_p \cap \cup_{\sigma \in \Sigma_p}\RR(\sigma_+^c),
\cr
\lambda(x) =0, 
\mbox{ on  } Y_p \setminus \cup_{\sigma \in \Sigma_p}\RR(\sigma_{++}^o), 
\cr
0 < \lambda(x) < 1, \mbox{ else}.  
\end{cases}
$$

We now construct a semi-algebraic homotopy, 
$$
h: Y_{p}  \times [0,\eps_{2p}] \rightarrow Y_p,
$$ 
by defining,

$$
\begin{array}{llll}
h(x,t) &=& \pi_{1\ldots k}\;\circ\;\phi(x, \lambda(x)t + 
(1 - \lambda(x))\varepsilon_{2p}),& \cr
&&\;\mbox{for}\; 0 < t \leq 
\eps_{2p},& \cr
h(x,0) &=& \lim_{t \rightarrow 0+} h(x,t), \;\mbox{else}.&
\end{array}
$$
Note that the last limit exists since $S$ is closed and bounded. 
We now show that, $h(x,0) \in \E(Z_{p+1},\R_{2p})$ for all $x \in Y_{p}.$

Let $x \in Y_p$ and $y = h(x,0)$. 

There are two  cases to consider.
\begin{itemize}
\item[$\lambda(x) < 1$:]
In this case, $x \in \E(Z_{p+1},\R_{2p})$ and by property (3) of $\phi$ and the
fact that $\lambda(x) < 1$, $y \in \E(Z_{p+1},\R_{2p})$.
\item[$\lambda(x) = 1$:]
Let $\sigma_y$ be the sign condition of ${\mathcal C}$ at $y$ and suppose that
$y \not\in \E(Z_{p+1},\R_{2p}).$
There are two cases to consider.
\begin{itemize}
\item[$\sigma_y \in \Sigma$:]
In this case, $ y \in X$ and hence there must exist 
$\tau \in{\rm Sign}_m({\mathcal C},S) \setminus \Sigma_m,$
with $m > p$ such that $y \in \RR(\tau_+^o).$
\item[$\sigma_y \not\in \Sigma$:]
In this case, taking $\tau = \sigma_y$, 
${\rm level}(\tau) > p$ and $y \in \RR(\tau_+^o).$
\end{itemize} 
It follows from the definition of $y,$ 
and property (3) of $\phi,$
that for any $m > p,$
and every $\rho \in {\rm Sign}_m({\mathcal C},S)$,
\begin{itemize}
\item
$y \in \RR(\rho_+^o)$ implies that $x \in \RR(\rho_+^o),$ and
\item
$x \in \RR(\rho_+^c)$ implies that $y \in \RR(\rho_+^c).$
\end{itemize}
Thus, $x \not\in Y_p$ which is a contradiction.
\end{itemize}

It follows that,
\begin{enumerate}
\item
$h(\cdot,\varepsilon_{2p}): Y_p \rightarrow Y_p$
is the identity map, 
\item
$h(Y_p,0) = \E(Z_{p+1},\R_{2p}),$ and
\item
$h(\cdot,t)$ restricted to $\E(Z_{p+1},\R_{2p})$ gives a semi-algebraic
homotopy between $h(\cdot,\varepsilon_{2p})|_{\E(Z_{p+1},\R_{2p})} = 
{\rm id}_{\E(Z_{p+1},\R_{2p})}$ and
$h(\cdot,0)|_{\E(Z_{p+1},\R_{2p})}$.
\end{enumerate}
Thus,
$Y_p$ is semi-algebraically homotopy equivalent to 
$\E(Z_{p+1},\R_{2p}).$
\end{proof}

\begin{lemma}
\label{lem:1}
For each $p, 0 \leq p \leq t,$
$Z_{p}$ is semi-algebraically homotopy equivalent to $\E(Y_{p},\R_{2p-1}).$  
\end{lemma}

\begin{proof}
For the purpose of the proof we define the following new sets
for $u \in \R_{2p}$.

\begin{enumerate}
\item
Let $Z_p'(u)\subset \R_{2p}^k$ be the set obtained by replacing 
in the definition of $Z_p$,
$\eps_{2j}$ by $\eps_{2j}-u$ and 
$\eps_{2j-1}$ by $\eps_{2j-1}+u$
for all $j>p,$ and 
$\eps_{2p}$ by $\eps_{2p}-u$, and 
$\eps_{2p-1}$ by $u$. 
For $u_0 > 0$ we will denote by $Z_p'((0,u_0])$ the set
$\{(x,u) \mid x \in Z_p'(u), u \in (0,u_0]\}.$
\item
Let $Y_p'(u) \subset \R_{2p}^k$ 
be the set obtained by replacing 
in the definition of $Y_p$,
$\eps_{2j}$ by $\eps_{2j}-u$ and
$\eps_{2j-1}$ by $\eps_{2j-1}+u$
for all $j > p$ and $\eps_{2p}$ by by $\eps_{2p}-u$.
\item
For $\sigma \in \Sign_m({\mathcal C},S)$, with $m \geq p$, let 
$\RR(\sigma_+^c)(u) \subset \R_{2p}^k$ 
denote  the set obtained by replacing 
$\eps_{2m}$ by $\eps_{2m}-u$
in the definition of $\RR(\sigma_+^c)$.
\item
For $\sigma \in \Sign_m({\mathcal C},S)$, with $m >  p$, let 
$\RR(\sigma_+^o)(u) \subset \R_{2p}^k$ 
denote  the set obtained by replacing 
$\eps_{2m-1}$ by $\eps_{2m-1}+u$
in the definition of $\RR(\sigma_+^o)$.
\item
Finally, for $\sigma \in \Sign_p({\mathcal C},S)$ let 
$\RR(\sigma_+^o)(u) \subset \R_{2p-1}^k$ denote 
the set obtained by replacing 
in the definition of $\RR(\sigma_o^c)$, 
$\eps_{2p-1}$ by $u$.
\end{enumerate}

Notice that by definition, for any $u,v \in \R_{2p}$ with  
$0 < u \leq v,$
$Z_p'(u) \subset Y_p'(u),$ 
$Z_p'(v) \subset Z_p'(u)$,
$Y_p'(v) \subset Y_p'(u)$, and
$$
\displaylines{
\bigcup_{0 < s \leq u} Y_p'(s) = \bigcup_{0 < s \leq u} Z_p'(s).
}
$$

We denote by $Z_p'$ (respectively, $Y_p'$) the
set $Z_p'(\eps_{2p-1})$ (respectively, $Y_p'(\eps_{2p-1})$).
It is easy to see that
$Y_p'$ is semi-algebraically homotopy equivalent to $\E(Y_p,\R_{2p-1})$,
and $Z_p'$ is semi-algebraically homotopy equivalent to $Z_p$.
We now prove that, $Y_p'$ is semi-algebraically
homotopy equivalent to $Z_p'$, which suffices to prove the lemma. 

Let 
$\mu: Y_{p}' \rightarrow \R_{2p-1}$ 
be the semi-algebraic map defined by 
\[
\mu(x) = \sup_{u \in (0,\eps_{2p-1}]} \{u \mid  x \in Z_{p}'(u) \}.
\]

We prove separately (Lemma \ref{lem:continuity} below)  
that $\mu$ is continuous. 
Note that the definition of the set $Z_p'(u)$ (as well as the set $Y_p'(u)$)
is more complicated than the more natural one consisting of just 
replacing $\eps_{2p-1}$ in the definition of $Z_p$ by $u$, 
is due to the fact that with
the latter definition the map $\mu$ defined below is not necessarily
continuous.

We now construct a continuous semi-algebraic map, 
$$
h: Y_p' \times [0,\varepsilon_{2p-1}] \rightarrow Y_p'
$$ 
as follows.

By Hardt's triviality theorem
there exist $u_0 \in \R_{2p},$ with $u_0 > 0$ and a semi-algebraic 
homeomorphism,
\[
\psi: Z_{p}'(u_0) \times (0,u_0] \rightarrow Z_p'((0,u_0]),
\]
such that 
\begin{enumerate}
\item
$\pi_{k+1}(\psi(x,u)) = u,$ 
\item
$\psi(x,u_0)  = (x,u_0)$ for $x \in Z'_{p}(u_0),$ 
and
\item
$\psi(\cdot,u)$ restricts to a homeomorphism of
$\RR(\sigma,Z_p'(u_0))$ to $\RR(\sigma,Z_p'(u))$, 
for every sign condition $\sigma$ of the
family, $\cup_{P \in {\mathcal C}} \{P, P\pm\eps_{2t},\ldots,P\pm
\varepsilon_{2p+1}\}$,  for all $u \in (0,u_0]$. 
\end{enumerate}

We now specialize $u_0$ to $\eps_{2p-1}$ and denote by $\phi$
the corresponding map,
\[
\phi: Z_{p}' \times (0,\eps_{2p-1}] \rightarrow Z_p'((0,\eps_{2p-1}]).
\]
Note, that for every $u$, $0 < u \leq \eps_{2p-1}$, $\phi$ gives
a homeomorphism,
$\phi_u: Z_p'(u) \rightarrow Z_p'$.
Hence, for every pair, $u,u'$,
$0 < u \leq u' \leq \eps_{2p-1}$, 
we have a homeomorphism,
$\theta_{u,u'} : Z_p'(u) \rightarrow \Z_p'(u')$ obtained by
composing $\phi_u$ with $\phi_{u'}^{-1}$. 
For $0 \leq u' < u \leq \eps_{2p-1}$, we let 
$\theta_{u,u'}$ be the identity map.
It is clear that $\theta_{u,u'}$ varies continuously with
$u$ and $u'$.

For $x \in Y_p', t \in [0,\varepsilon_{2p-1}]$ we now define,
\[
h(x,t) = \theta_{\mu(x),t}(x).
\]

It is easy to verify from the definition of $h$ 
and the properties of $\phi$ listed above that, 
$h$ is continuous and satisfies the following.
\begin{enumerate}
\item
$h(\cdot,0): Y_p' \rightarrow 
Y_p'$ is the identity map, 
\item
$h(Y_p',\varepsilon_{2p-1}) = Z_{p}',$ and
\item
$h(\cdot,t)$ restricts to a homeomorphism 
$Z_{p}' \times t \rightarrow Z_{p}'$
for every $t \in [0,\varepsilon_{2p-1}].$
\end{enumerate}
This proves the required homotopy equivalence.
\end{proof}

We now prove that the function $\mu$ used in the proof above is continuous. 

\begin{lemma}
\label{lem:continuity}
The semi-algebraic map
$\mu: Y_{p}' \rightarrow \R_{2p-1}$ 
defined by 
\[
\mu(x) = \sup_{u \in (0,\eps_{2p-1}]} \{u \mid  x \in Z_{p}'(u) \}
\]
is continuous.

\end{lemma}

\begin{proof}
Let $0 < \delta \ll \eps_{2p-1}$ be a new infinitesimal. 
In order to prove the continuity of $\mu$
(which is a semi-algebraic function defined over $\R_{2p-1}$),
it suffices to show that 
$$
\lim_\delta \E(\mu, \R_{2p-1}\la\delta\ra)(x') = 
\lim_\delta \E(\mu, \R_{2p-1}\la\delta\ra)(x)
$$
for every pair of points
$x,x' \in \E(Y_p',\R_{2p-1}\la\delta\ra)$ 
such that $\lim_\delta x = \lim_\delta x'$.

Consider such a pair of points $x,x' \in \E(Y_p',\R_{2p-1}\la\delta\ra).$ 
Let $u \in (0,\eps_{2p-1}]$ be such that
$x \in Z_p'(u)$.
We show below  that this implies $x' \in Z_p'(u')$ for
some $u'$ satisfying $\lim_\delta u' = \lim_\delta u.$

Let $m$ be the largest integer such that
there exists $\sigma \in \Sigma_m$ with 
$x \in \RR(\sigma_+^c)(u)$.
Since $x \in Z_p'(u)$ such an $m$ must exist.

We have two cases:
\begin{enumerate}
\item
$m > p$: Let $\sigma \in \Sigma_m$ with 
$x \in \RR(\sigma_+^c)(u).$ Then, by the maximality of $m$,
we have that for each $P \in {\mathcal C}$,
$\sigma(P) \neq 0$ implies that $\lim_\delta P(x) \neq 0$. 
As a result, we have that $x' \in \RR(\sigma_+^c)(u')$ 
for all $u' < u - \max_{P \in {\mathcal P}, \sigma(P) = 0} |P(x) - P(x')|,$ 
and hence we can choose $u'$ such that 
$x' \in \RR(\sigma_+^c)(u')$  and  $\lim_\delta u' = \lim_\delta u$.
\item
$m \leq p$:
If  $x' \not\in Z_p'(u)$ then since $x' \in Y_p' \subset Y_p'(u)$, 
$$
x' \in \cup_{\sigma \in {\rm Sign}_p({\mathcal C},S) \setminus \Sigma_p}
\RR(\sigma_+^o)(u).
$$  
Let $\sigma \in {\rm Sign}_p({\mathcal C},S) \setminus \Sigma_p$
be such that $x' \in \RR(\sigma_+^o)(u)$.
We prove by contradiction that 
$\displaystyle{
\lim_\delta \max_{P \in {\mathcal P},\sigma(P) = 0}|P(x')| = u.
}$

Assume that 
$$\lim_\delta \max_{P \in {\mathcal P},\sigma(P) = 0}|P(x')| \neq u.$$
Since, $ x \not\in  \RR(\sigma_+^o)(u)$ by assumption, and
$\lim_\delta x' = \lim_\delta x$, 
there must exist $P \in {\mathcal C}$, $\sigma(P) \neq 0$, and 
$\lim_\delta P(x) = 0$.
Letting $\tau$ denote the sign condition defined by 
$\tau(P) = 0$ if  $\lim_\delta P(x) = 0$ and
$\tau(P) = \sigma(P)$ else,
we have that ${\rm level}(\tau) > p$ and $x$ belongs to both
$\RR(\tau_+^o)(u)$ as well as $\RR(\tau_+^c)(u)$.

Now there are two cases to consider depending on whether
$\tau$ is in $\Sigma$ or not. If $\tau \in \Sigma$, then
the fact that $x \in \RR(\tau_+^c)(u)$ contradicts the choice of 
$m$, since $m \leq p$ and ${\rm level}(\tau) > p$.
If $\tau \not\in\Sigma$ then $x$ gets removed at the level
of $\tau$ in the construction of $Z_p'(u)$, and hence 
$x \in \RR(\rho_+^c)(u)$ for some $\rho \in \Sigma$
with ${\rm level}(\rho) > {\rm level}(\tau) > p$. This again
contradicts the choice of $m$. 
Thus,
$\displaystyle{
\lim_\delta \max_{P \in {\mathcal P},\sigma(P) = 0} |P(x')| = u
}$ 
and since
$\displaystyle{
x' \not\in \cup_{\sigma \in {\rm Sign}_p({\mathcal C},S) \setminus \Sigma_p}
\RR(\sigma_+^o)(u')
}$ 
for all 
$\displaystyle{
u' < \max_{P \in {\mathcal P},\sigma(P) = 0}|P(x')|,
}$
we can choose 
$u'$ such that 
$\lim_\delta u' = \lim_\delta u$, and
$
\displaystyle{
x' \not\in \cup_{\sigma \in {\rm Sign}_p({\mathcal C},S) \setminus \Sigma_p}
\RR(\sigma_+^o)(u').
}$
\end{enumerate}
In both cases we have that  $x' \in Z_p'(u')$ for
some $u'$ satisfying $\lim_\delta u' = \lim_\delta u,$
showing that $\lim_\delta\mu(x') \geq \lim_\delta \mu(x).$
The reverse inequality follows by exchanging the roles of $x$ and $x'$ in
the previous argument.
Hence, $\lim_\delta \mu(x') = \lim_\delta \mu(x),$
proving the continuity of $\mu$.
\end{proof}

\begin{proof}(Proof of Theorem \ref{the:GV})
The theorem follows immediately from Lemmas \ref{lem:2} and \ref{lem:1}.
\end{proof}

\section{Computing the first Betti number of a general semi-algebraic set}
\label{sec:general}
In this section we describe the algorithm for computing the first 
Betti number of a general semi-algebraic set. We first replace the
given set by a closed and bounded one, using the construction described
in the previous section. We then apply
Algorithm \ref{alg:firstbetticlosed}.

\begin{algorithm}
(First Betti Number of a ${\mathcal P}$- Semi-algebraic Set)
\index{Connected components!Algorithm!Algebraic Set}
\begin{itemize}
\item[]
\item[] {\sc Input.} 
\begin{itemize}
\item  
a finite set of  polynomials ${\mathcal P} \subset \D[X_1,\ldots,X_k],$
\item a formula defining a ${\mathcal P}$-semi-algebraic set, $S$.
\end{itemize}

\item[]{\sc Output.}
$b_1(T).$
\item[] {\sc Procedure.}
\item[Step 1] Let $\eps$ be an infinitesimal. Define $\tilde S$ as the intersection of
$\E(S,\la\eps\ra)$ with the ball of center $0$ and radius $1/\eps$.
Define ${\mathcal Q}$ as ${\mathcal P} \cup \{\eps^2(X_1^2+\ldots+X_k^2+X_{k+1}^2)-4, X_{k+1}\}$ 
 Replace $\tilde S$ by the ${\mathcal Q}$- semi-algebraic set $S$ defined as the intersection
 of the cylinder $\tilde S\times  \R\la \eps \ra$ with the upper hemisphere defined by
 $\eps^2(X_1^2+\ldots+X_k^2+X_{k+1}^2)=4, X_{k+1}\ge 0.$
\item[Step 2]
Using the Gabrielov-Vorobjov construction described above,   
replace $T$ by   a ${\mathcal Q}'$-closed set, $T'$.
\item[Step 3]
Use Algorithm \ref{alg:firstbetticlosed}
to compute the first Betti number of $T'$.
\item[] {\sc Complexity.}
The complexity of the algorithm is bounded by 
$(sd)^{k^{O(1)}}$,
where $s = \#{\mathcal P}$ and $d = \max_{P \in {\mathcal P}}\deg(P).$
\end{itemize}
\end{algorithm}

\vspace{.1in}
\noindent
{\sc Proof of correctness.}
The correctness of the algorithm is a consequence of 
Theorem \ref{the:GV} and the correctness of
Algorithm \ref{alg:firstbetticlosed}.
\eop

\vspace{.1in}
\noindent
{\sc Complexity analysis.}
In Step 2 of the algorithm the cardinality of ${\mathcal Q}'$ is 
$2(s+1)^2$ and the degrees of the
polynomials in ${\mathcal Q}'$ are still bounded by $d$.
The complexity of Step 3 of the algorithm is then bounded by
$(sd)^{k^{O(1)}}$ using the complexity analysis of  
Algorithm \ref{alg:firstbetticlosed}.
\eop

\section{Computing Connected Components}
\label{sec:cc}
If one is  
interested in computing semi-algebraic descriptions of the
connected components of a given semi-algebraic set, then 
using Algorithm \ref{16:alg:paramconnecting} (Parametrized Bounded Connecting)
it is possible to do so with a complexity 
making precise the one of
previously known algorithms, whose complexities were of the form
$(sd)^{k^{O(1)}}$ (see \cite{HRS94}).
We have the following theorems (we refer the reader 
to \cite{BPR03} for details of the proof).

\begin{theorem}
\label{15:the:algebraiccc}
 If $\ZZ(Q,\R^k)$ is an algebraic set
 defined as the zero set of a polynomial $Q \in
{\D}[X_1,\ldots,X_k]$ of degree $\leq d$, then there
is an algorithm that outputs
quantifier free formulas whose realizations are the
 semi-algebraically connected
 components of $\ZZ(Q,\R^k).$ The complexity of
the algorithm
in the ring generated by the coefficients of $Q$
is bounded by $ d^{O(k^3)}$ and the
 degrees of the polynomials that appear in the
output are bounded by
$O(d)^{k^2}.$
Moreover, if $\D=\Z,$ and the bitsizes
of the coefficients of the polynomials are bounded by
 $\tau$, then the bitsizes of the integers appearing in the
intermediate computations and the output are bounded
by $\tau  d^{O(k^2)}$.
\end{theorem}

\begin{theorem}
\label{16:the:dcc}
Let
${\mathcal P} =
\{P_1,\ldots,P_s \} \subset {\D}[X_1,\ldots,X_k]$
 with $\deg(P_i) \leq d,
1 \leq i \leq s$
and a semi-algebraic set $S$
 defined by a ${\mathcal
P}$
quantifier-free  formula.
There exists  an algorithm that outputs
quantifier-free
 semi-algebraic descriptions of all the
semi-algebraically connected components of
$S$.
The
complexity of the algorithm is bounded by
$s^{k+1} d^{O(k^4)} .$
The degrees of the polynomials that
appear in the output are bounded by $d^{O(k^3)} .$
Moreover, if the input polynomials have integer
coefficients whose bitsize is bounded by $\tau$ the bitsize of coefficients
output is
$d^{O(k^3)}\tau$.
\end{theorem}


\begin{thebibliography}{50}
\bibitem{B99}
{\sc S.\ Basu},
\newblock {\em On Bounding the Betti Numbers and Computing the Euler
Characteristics of Semi-algebraic Sets},
\newblock {Discrete and Computational Geometry}, 22:1-18 (1999).

\bibitem{B03}
{\sc S.\ Basu},
\newblock{\em On different bounds on different Betti numbers},
\newblock {Discrete and Computational Geometry}, 
30:1, 65-85, (2003).

\bibitem{B04}
{\sc S.\ Basu},
\newblock{\em Computing the first few Betti numbers of semi-algebraic sets
in singly exponential time},
\newblock preprint (2004).
(Available at {\tt www.math.gatech.edu/$\sim$saugata/bettifew.ps}.)

\bibitem{BPR95}
{\sc S.\ Basu, R.\ Pollack, M.-F. \ Roy},
 \newblock {\em On the
Combinatorial and Algebraic
Complexity of Quantifier Elimination},
\newblock Journal of the ACM , 43  1002--1045, (1996).

\bibitem{BPR99}
{\sc S.\ Basu, R.\ Pollack, M.-F. \ Roy},
 \newblock {\em Computing Roadmaps of Semi-algebraic
Sets on a Variety},
\newblock Journal of the AMS, vol 3, 1 55-82 (1999).

\bibitem{BPR03}
{\sc S.\ Basu, R.\ Pollack, M.-F. \ Roy},
\newblock {\em Algorithms in Real Algebraic Geometry},
Springer-Verlag, 2003. 
(Updated version available electronically at: \\
{\tt www.math.gatech.edu/$\sim$saugata/bpr-posted1.html}.)

\bibitem{BCR}
{\sc J.\ Bochnak, M.\ Coste, M.-F.\ Roy},
\newblock {\em G\'eom\'etrie
alg\'ebrique r\'eelle,}
Springer-Verlag (1987).
\newblock {\em Real algebraic
geometry},
Springer-Verlag (1998).

\bibitem{BM}
{\sc A. \ Borel, J.C.\ Moore},
\newblock {\em Homology theory for locally compact spaces},
\newblock {Mich. Math. J.}, 7:137-- 159, (1960).

\bibitem{BC}
{\sc P.\ Burgisser, F.\ Cucker},
\newblock {\em Counting Complexity Classes for Numeric Computations II: 
Algebraic and  Semi-algebraic Sets},
\newblock preprint.

\bibitem{Canny93a}
{\sc J.\ Canny},
\newblock {\em Computing road maps in general semi-algebraic sets},
\newblock The
Computer Journal, 36: 504--514, (1993).

\bibitem{Col}
{\sc G. Collins},
\newblock {\em Quantifier elimination for real closed fields by
cylindric algebraic decomposition},
\newblock In
Second GI Conference on Automata Theory and
Formal Languages. Lecture
Notes in Computer Science, vol. 33, pp. 134-183, Springer-
Verlag, Berlin (1975).

\bibitem{GV}
{\sc A.\ Gabrielov, N.\ Vorobjov}
\newblock {Betti Numbers for Quantifier-free Formulae,}
\newblock Discrete and Computational Geometry, 33:395-401, 2005.

\bibitem{GR92}
{\sc L. Gournay, J. J. Risler},
\newblock{\em Construction of roadmaps of semi-algebraic sets},
\newblock
Appl. Algebra Eng. Commun. Comput. 4, No.4, 239-252 (1993).


\bibitem{GV92}
{\sc D. Grigor'ev, N. Vorobjov},
\newblock{\em Counting connected components of a semi-algebraic
set in sub-exponential time},
\newblock Comput. Complexity 2, No.2, 133-186 (1992).


\bibitem{Hardt} {\sc R. M. Hardt},
\newblock{\em Semi-algebraic Local Triviality in Semi-algebraic Mappings},
\newblock Am. J. Math. { 102}, 291-302 (1980).

\bibitem{HRS94}
{\sc J.\ Heintz, M.-F.\ Roy, P. Solern\`o},
\newblock{\em Description of the
Connected Components of a Semi-algebraic Set in Single Exponential Time},
\newblock{Discrete and Computational Geometry}, 11, 121-140 (1994).

\bibitem{Milnor}
{\sc J. \ Milnor},
\newblock {\em On the Betti numbers of real varieties},
\newblock Proc. AMS 15, 275-280 (1964).

\bibitem{O}
{\sc O.\ A.\ Ole\u{i}nik},
\newblock {\em Estimates of the {B}etti numbers
of real algebraic hypersurfaces},
\newblock { Mat. Sb. (N.S.)}, 28 (70): 635--640 (Russian) (1951).

\bibitem{OP}
{\sc O. A.\ Ole\u{i}nik, I. B.\ Petrovskii},
\newblock {\em On the topology of real algebraic surfaces},
\newblock Izv. Akad. Nauk SSSR 13, 389-402 (1949).

\bibitem{R92}
{\sc J.\ Renegar.}
\newblock {\em On the computational complexity and geometry of the
first order theory of the reals},
\newblock Journal  of Symbolic Computation, 13: 255--352 (1992).

\bibitem{Thom}
{\sc R.\ Thom},
\newblock {\em Sur l'homologie des vari\'et\'es alg\'ebriques
r\'eelles},
\newblock  {Differential and Combinatorial
Topology},  255--265.
Princeton University Press, Princeton (1965).
\end{thebibliography}
\end{document}